\documentclass[hidelinks,onefignum,onetabnum]{siamart250211}
\usepackage{aligned-overset}
\usepackage{amsmath,amsfonts,amssymb}
\usepackage[style=numeric-comp,maxbibnames=99,bibencoding=utf8,giveninits=true]{biblatex}
\usepackage{booktabs}
\usepackage{cancel}
\usepackage{hyperref}
\usepackage{xcolor}

\usepackage{pgfplots}
\usepackage[outline]{contour}
\usepackage{subcaption}
\usepackage{multirow}

\bibliography{hypre-rsr}

\ifpdf
\hypersetup{
  pdftitle={A Reynolds-semi-robust method with hybrid velocity and pressure for the unsteady incompressible Navier--Stokes equations},
  pdfauthor={L. Beir\~{a}o da Veiga, D. A. Di Pietro, J. Droniou, K. B. Haile, T. J. Radley}
}
\fi

\headers{A Reynolds-semi-robust method with hybrid velocity and pressure}{L. Beir\~{a}o da Veiga, D. A. Di Pietro, J. Droniou, K. B. Haile, T. J. Radley}

\newsiamremark{remark}{Remark}

\newcommand{\st}{\;:\;}

\newcommand{\Real}{\mathbb{R}}

\newcommand{\Th}{\mathcal{T}_h}

\newcommand{\Fh}{\mathcal{F}_h}
\newcommand{\Fhb}{\mathcal{F}_h^{\rm b}}
\newcommand{\FT}{\mathcal{F}_T}

\newcommand{\Poly}[1]{\mathcal{P}^{#1}}
\newcommand{\RTN}[1]{\mathcal{RT\!N}^{#1}}

\newcommand{\UT}{\underline{U}_T^k}
\newcommand{\Uh}{\underline{U}_h^k}
\newcommand{\UhZ}{\underline{U}_{h,0}^k}
\newcommand{\ZhZ}{\underline{Z}_{h,0}^k}

\newcommand{\PT}{\underline{P}_T^k}
\newcommand{\Ph}{\underline{P}_h^k}
\newcommand{\PhZ}{\underline{P}_{h,0}^k}

\newcommand{\IUh}{\underline{I}_{U,h}^k}
\newcommand{\IPh}{\underline{I}_{P,h}^k}
\newcommand{\IUT}{\underline{I}_{U,T}^k}
\newcommand{\IPT}{\underline{I}_{P,T}^k}
\newcommand{\IRTNT}[1]{I_{\mathcal{RT\!N},T}^{#1}}

\newcommand{\lproj}[2]{\pi^{#1}_{#2}}

\newcommand{\GT}{G_T^k}
\newcommand{\RT}{R_T^{k+1}}
\newcommand{\RTo}{R_T^1}
\newcommand{\Rh}{R_h^{k+1}}

\newcommand{\tF}{t_{\rm F}}

\newcommand{\Hdiv}[1]{H(\operatorname{div};#1)}

\newcommand{\norm}[2]{\|#2\|_{#1}}
\newcommand{\seminorm}[2]{|#2|_{#1}}
\newcommand{\Norm}[2]{\left\|#2\right\|_{#1}}
\newcommand{\Seminorm}[2]{\left|#2\right|_{#1}}
\newcommand{\tnorm}[2]{|\kern-0.25ex|\kern-0.25ex|#2|\kern-0.25ex|\kern-0.25ex|_{#1}}

\newcommand{\term}{\mathfrak{T}}

\newcommand{\huline}[1]{\widehat{\underline{#1}}}

\newcommand{\Err}{\mathcal{E}}

\newcommand{\ReT}{\mathrm{Re}_T}

\begin{document}

\title{A Reynolds-semi-robust method with hybrid velocity and pressure for the unsteady incompressible Navier--Stokes equations}
\author{L. Beir\~{a}o da Veiga\thanks{Dipartimento di Matematica e Applicazioni, Università degli Studi di Milano-Bicocca,
    Piazza dell’Ateneo Nuovo 1, 20126 Milano, Italy (\email{lourenco.beirao@unimib.it})}
  \and D. A. Di Pietro\thanks{IMAG, Univ Montpellier, CNRS, Montpellier 34090, France (\email{daniele.di-pietro@umontpellier.fr})},
  \and J. Droniou\thanks{IMAG, Univ Montpellier, CNRS, Montpellier 34090, France (\email{jerome.droniou@umontpellier.fr}), School of Mathematics, Monash University, Melbourne, Australia},
  \and K. B. Haile\thanks{Dipartimento di Matematica e Applicazioni, Università degli Studi di Milano-Bicocca,
    Piazza dell’Ateneo Nuovo 1, 20126 Milano, Italy (\email{k.haile@campus.unimib.it})}
  \and T. J. Radley \thanks{IMAG, Univ Montpellier, CNRS, Montpellier 34090, France (\email{thomas.radley1@umontpellier.fr})}
}

\maketitle

\begin{abstract}
  In this paper we propose and analyze a new Finite Element method for the solution of the two- and three-dimensional incompressible Navier--Stokes equations based on a hybrid discretization of both the velocity and pressure variables.
  The proposed method is pressure-robust, i.e., irrotational forcing terms do not affect the approximation of the velocity, and Reynolds-quasi-robust, with error estimates that, for smooth enough exact solutions, do not depend on the inverse of the viscosity.
  We carry out an in-depth convergence analysis highlighting pre-asymptotic convergence rates and validate the theoretical findings with a complete set of numerical experiments.
\end{abstract}



\section{Introduction}

The development and analysis of Finite Element Methods (FEM) for the incompressible Navier--Stokes equations (both stationary and time-dependent) have long been a central focus of research in the mathematical community. Alongside the specific articles referenced below, we highlight as examples the monographs \cite{Girault.Raviart:86,quartapelle:book,volker:book} or the article \cite{VKN18} providing modern perspectives.

Among the numerous contributions, two key numerical challenges in the modern FEM literature on the Navier--Stokes equations stand out. The first, and historically the most significant, is the difficulty posed by convection-dominated flows, which can result in poor convergence of numerical schemes and non-physical oscillations in the discrete solution. To address cases where convection prevails over diffusion, various stabilization techniques have been proposed. Limiting our citations to a few notable references, we mention the well-known Streamline Upwind Petrov--Galerkin (SUPG) method and its variants \cite{FF:1982,BH:1982,TB:1996,Beirao-da-Veiga.Dassi.ea:23}, the Continuous Interior Penalty (CIP) method \cite{BFH:2006,BF:2007}, Grad-Div stabilization \cite{OLHL:2009,grad-div}, Local Projection (LP) \cite{BB:2006,MT:2015,LPS-NS}, and the Variational MultiScale approach \cite{JK:2005,HJ:2000}. In this respect, we refer to a scheme as \emph{convection quasi-robust} if, assuming that the exact solution is sufficiently regular, the velocity error is independent of inverse powers of the diffusion coefficient, possibly in a norm which includes also some control on convection.

A second, more recent, aspect of importance is {\it pressure robustness} \cite{Linke:14}. In essence, a pressure-robust scheme ensures that modifications to the continuous problem that affect only the pressure result in changes to the discrete pressure only, leaving the discrete velocity unaltered. This property, thoroughly investigated in several studies, offers various advantages, as discussed in \cite{Linke.Merdon:16,John.Linke.ea:17}. One way to achieve pressure robustness is to employ a FEM that guarantees a divergence-free discrete kernel, assuming the forms involved are approximated exactly, or possibly by using specialized projections \cite{Linke:14}. We also point to the discussions on Virtual Element methods in \cite{Beirao-da-Veiga.Lovadina.ea:18,BMV:2018} and on Hybrid-High Order methods in \cite{Di-Pietro.Ern.ea:16*1,Castanon-Quiroz.Di-Pietro:20,Castanon-Quiroz.Di-Pietro:23}; see also \cite{Beirao-da-Veiga.Dassi.ea:22,Di-Pietro.Droniou.ea:24} for an approach inspired by a non-standard weak formulation.

Of particular relevance for the present article are Discontinuous Galerkin (DG) methods which, as it first appeared evident for the simpler advection-diffusion problem \cite{Di-Pietro.Ern.ea:08,Ayuso.Marini:09,Ern.Stephansen.ea:09}, are particularly well-suited for this kind of equations due to the higher convection control yielded by discontinuous test functions. Indeed, using DG schemes allows for a very simple and effective stabilization based on upwinding. In the framework of incompressible fluids, ideas stemming from the DG methodology allowed to make use of finite elements borrowed from the de Rham complex \cite{Hdiv1,Hdiv2,Hdiv3}; due to the exact enforcement of the zero-divergence condition, such construction resulted in novel schemes which are naturally pressure robust. Furthermore, thanks to the non $H^1$-conforming nature of the discrete spaces, these methods can also be stabilized through upwinding, leading to a class of schemes which are currently considered among the most successful FEMs for incompressible fluids at the academic level.

In the present work we propose a novel FEM for the discretization of the time-dependent incompressible Navier--Stokes equation in $2$ and $3$ space dimensions extending the construction for the steady Stokes problem proposed in \cite[Section~5.4]{Botti.Botti.ea:24} (see also \cite{Botti.Massa:22} for a related method).
Specifically, this method is recovered from the present one dropping the unsteady and convective terms and selecting $\Sigma_T = \nabla \mathcal{P}^{k+1}(T)^d$ in the notation of the above reference.
A careful design of the spaces and of the convective term makes the scheme both convection quasi-robust and pressure-robust.
Both the velocity and the pressure have a set of internal degrees of freedom, associated to elements, and a set of boundary degrees of freedom, associated to faces.
In this sense, the proposed approach has similarities with hybridized DG schemes for Navier--Stokes using de Rham complex spaces (see Remark \ref{rem:keegan.et.al}).
By a careful construction of the discrete convection form, non-dissipativity is acquired without the need of an artificial anti-symmetrization.
Furthermore, by making use of a high-order reconstruction in the spirit of Hybrid High-Order methods and Virtual Elements Methods, we are able to develop a more accurate approximation of the diffusive term with respect to the natural polynomial order of the original spaces.
Our construction is combined with an upwind-like stabilization, that here takes the form of a penalization of the jumps among the face velocity variables and the trace of the internal velocity variables.
Through the use of regime-dependent estimates of the convective component of the consistency error in the spirit of \cite{Di-Pietro.Droniou.ea:15,Botti.Di-Pietro.ea:18,Di-Pietro.Droniou:23*2,Beirao-da-Veiga.Di-Pietro.ea:24}, our $h$-convergence analysis accounts for pre-asymptotic orders of convergence in convection-dominated configurations; see Corollary~\ref{cor:convergence.rate} and Remark~\ref{rem:convergence.rate} concerning the diffusion-dominated case.

An important asset in the proposed approach is the potential for efficient extensions to non-linear constitutive laws.
Indeed, in order to be extended to more general (nonlinear) rheological fluid laws, DG schemes typically need to make use of suitable gradient reconstruction operators, which in turn lead to very large stencils in the involved matrices; see, e.g., \cite{Burman.Ern:08,Del-Pezzo.Lombardi.ea:12,Diening.Koner.ea:14,Malkmus.Ruzicka.ea:18,Beirao-da-Veiga.Di-Pietro.ea:24} concerning $p$-type diffusion.
Thanks to the use of face variables, the present approach paves the way for schemes with much smaller stencils, still preserving all the main advantages of DG FEMs.

The rest of the article is organized as follows. After presenting the continuous problem and some preliminaries in Section \ref{sec:prelimins}, in Section \ref{sec:scheme} we describe the discrete scheme and show its well-posedness. In Section \ref{sec:error-analysis} we develop the convergence analysis of the method, proving a convection quasi-robust error estimate for the velocity variable in a natural energy-like discrete norm.
Furthermore, the error bound is independent of the pressure variable, resulting in pressure robustness of the scheme. Finally, a set of numerical tests are presented in Section~\ref{sec:num} to corroborate the theoretical results.


\section{Notation and preliminary results}\label{sec:prelimins}

In the present section we briefly review the continuous problem and introduce some preliminary definitions and results which will be instrumental to the following developments.

\subsection{Continuous problem}

Let $\Omega \in \Real^d$, $d \in \{2, 3\}$, denote an open bounded connected polygonal (if $d = 2$) or polyhedral (if $d = 3$) domain with Lipschitz boundary $\partial \Omega$ and outward unit normal vector $n$, and let $t_F>0$ be a set final time. The unsteady incompressible Navier--Stokes problem reads as follows:
Given $f : (0,\tF\rbrack \times \Omega \to \Real^d$ and $u_0 : \Omega \to \Real^d$ such that $u_0 = 0$ on $\partial \Omega$ and $\nabla \cdot u_0 = 0$, find the velocity $u : [0,\tF] \times \Omega \to \Real^d$ and pressure $p : (0,\tF\rbrack \times \Omega \to \Real$ satisfying
\begin{subequations}\label{eq:strong}
  \begin{equation}\label{eq:strong:ic}
    u(0,\cdot) = u_0
  \end{equation}
  and, for $t \in (0,\tF\rbrack$,
  \begin{alignat}{2}\label{eq:strong:momentum}
    \partial_t u(t) - \nu \Delta u(t) + (u(t) \cdot \nabla) u(t) + \nabla p(t) &= f(t) &\qquad& \text{in $\Omega$},
    \\ \label{eq:strong:mass}
    \nabla \cdot u(t) &= 0 &\qquad& \text{in $\Omega$},
    \\ \label{eq:strong:bc}
    u(t) &= 0 &\qquad& \text{on $\partial \Omega$},
    \\ \label{eq:strong:p.zero-average}
    \int_\Omega p(t) &= 0,
  \end{alignat}
\end{subequations}
where $\nu > 0$ is the kinematic viscosity and, given a function of time and space $\psi$, we have adopted the convention that $\psi(t)$ stands for the function of space $\psi(t,\cdot)$.

\subsection{Mesh}

Let $\Th$ be a matching simplicial mesh of $\Omega$ belonging to a regular sequence in the sense of \cite{Ciarlet:02}, and denote by $\Fh$ the corresponding set of simplicial faces.
Boundary faces contained in $\partial \Omega$ are collected in the set $\Fhb$.
Given a mesh element $T \in \Th$, we denote by $\FT$ the subset of $\Fh$ collecting the faces that lie on its boundary and, for all $F \in \FT$, $n_{TF}$ is the unit vector normal to $F$ and pointing out of $T$.

The diameter of a mesh element or face $Y \in \Th \cup \Fh$ is denoted by $h_Y$, so that $h = \max_{T \in \Th} h_T$.
To avoid the proliferation of generic constants, from this point on we abbreviate as $a \lesssim b$ the inequality $a \le C b$ with real number $C > 0$ independent of the meshsize $h$, of the viscosity $\nu$ and, for local inequalities on a mesh element or face $Y \in \Th \cup \Fh$, of $Y$, but possibly depending on other parameters such as the mesh regularity parameter, the polynomial degree, the domain, and $t_F$. When relevant, the dependencies of $C$ will be specified more precisely.
We also write $a \simeq b$ for ``$a \lesssim b$ and $b \lesssim a$''.


\subsection{Polynomial spaces}

Given a mesh element or face $Y \in \Th \cup \Fh$ and an integer $\ell \ge 0$, we denote by $\Poly{\ell}(Y)$ the space spanned by the restriction to $Y$ of polynomial functions of the space variables of total degree $\le \ell$ and we set, by convention, $\Poly{-1}(Y) \coloneq \{0\}$.
The $L^2$-orthogonal projector on $\Poly{\ell}(Y)$ is $\lproj{\ell}{Y} : L^1(Y) \to \Poly{\ell}(Y)$ such that, for all $q \in L^1(Y)$,
\begin{equation}\label{eq:lproj}
  \int_Y \lproj{\ell}{Y} q \, r = \int_Y q \, r
  \qquad \forall r \in \Poly{\ell}(Y).
\end{equation}
When applied to vector-valued functions, $\lproj{\ell}{Y}$ acts component-wise.

Let now $T \in \Th$.
The Raviart--Thomas--Nédélec space of degree $\ell \ge 1$ is
\[
\RTN{\ell}(T) \coloneq \Poly{\ell-1}(T)^d + x \Poly{\ell-1}(T),
\]
with interpolator $\IRTNT{\ell} : H^1(T)^d \to \RTN{\ell}(T)$ such that, for all $v \in H^1(T)^d$,
\begin{equation}\label{eq:IRTNT}
  \text{
    $\lproj{\ell-2}{T} \IRTNT{\ell} v = \lproj{\ell-2}{T} v$
    \quad and \quad $\IRTNT{\ell} v \cdot n_{TF} = \lproj{\ell-1}{F} ( v \cdot n_{TF} )$ for all $F \in \FT$.
  }
\end{equation}
We note the following key commutation property:
\begin{equation}\label{eq:IRTNT:commutation}
  \nabla \cdot (\IRTNT{\ell} v) = \lproj{\ell-1}{T} (\nabla \cdot v)
  \qquad \forall v \in H^1(T)^d.
\end{equation}

\begin{lemma}[Approximation properties of the $\RTN{\ell}$ interpolator]\label{lem:approx.RTN}
  For all $s\in [1,\infty]$, all $\ell\ge 1$, and all integers $0 \le q \le \ell-1$ and $0 \le m \le q + 1$, it holds
  \begin{equation}\label{eq:IRTNT:approximation}
    \seminorm{W^{m,s}(T)^d}{v - \IRTNT{\ell} v}
    \lesssim h_T^{q+1-m} \seminorm{W^{q+1,s}(T)^d}{v}
    \qquad \forall v \in W^{q+1,s}(T)^d
  \end{equation}
  and, for all $F\in\FT$ and $\alpha\in\mathbb{N}^d$ with $\sum_{i=1}^d\alpha_i \eqcolon r\le q$,
  \begin{equation}\label{eq:IRTNT:approximation.trace}
    \norm{L^s(F)^d}{\partial^\alpha(v - \IRTNT{\ell} v)}
    \lesssim h_T^{q+1-r-\frac1s} \seminorm{W^{q+1,s}(T)^d}{v}
    \qquad \forall v \in W^{q+1,s}(T)^d.
  \end{equation}
\end{lemma}

\begin{remark}[Boundedness of the $\RTN{\ell}$ interpolator]
  A triangle inequality and \eqref{eq:IRTNT:approximation} with $m=q+1$ immediately yield, for the same range of indices $(s,\ell,q)$ as in Lemma \ref{lem:approx.RTN},
  \begin{equation}\label{eq:boundedness.IRTNT}
    \seminorm{W^{q+1,s}(T)^d}{\IRTNT{\ell}v}\lesssim \seminorm{W^{q+1,s}(T)^d}{v} .
  \end{equation}
\end{remark}

\begin{proof}[Proof of Lemma \ref{lem:approx.RTN}]
  In the Hilbertian case $s=2$, this result is proved, e.g., in \cite[Lemma~3.17]{Gatica:14};
  see also \cite[Proposition~2.5.1]{Boffi.Brezzi.ea:13}, or \cite[Theorem~16.4]{Ern.Guermond:21*1}.
  As we will require it for both $s = 2$ and $s = \infty$, we provide an independent proof in generic Sobolev spaces.
  An argument based on the use of a reference element (through the Piola transform \cite[Section 3.4.2]{Gatica:14}) and the unisolvence of the degrees of freedom of the Raviart--Thomas--N\'edelec element \cite[Theorem 3.3]{Gatica:14} easily gives
  \[
  \norm{L^s(T)^d}{z}\lesssim \norm{L^s(T)^d}{\lproj{\ell-2}{T}z}+h_T^{\frac1s}\sum_{F\in\FT}\norm{L^s(F)}{z\cdot n_{TF}}\qquad\forall z\in\RTN{\ell}(T).
  \]
  Taking $w\in W^{1,s}(T)$, applying this estimate to $z = \IRTNT{\ell}w$, and recalling \eqref{eq:IRTNT}, we obtain
  \begin{equation}\label{eq:IRTNT:boundedness}
    \begin{aligned}
      \norm{L^s(T)^d}{\IRTNT{\ell}w}
      \lesssim{}&
      \norm{L^s(T)^d}{\lproj{\ell-2}{T}w}
      + h_T^{\frac1s}\sum_{F\in\FT}\norm{L^s(F)}{\lproj{\ell-1}{F}(w\cdot n_{TF})}
      \\
      \lesssim{}&
      \norm{L^s(T)^d}{w}+h_T\norm{L^s(T)^{d\times d}}{\nabla w},
    \end{aligned}
  \end{equation}
  where we have used the $L^s$-boundedness of $\lproj{\ell-2}{T}$ and $\lproj{\ell-1}{F}$ followed by a local trace inequality (see, e.g., \cite[Lemma~12.8]{Ern.Guermond:21*1} and also \cite[Lemma 1.44 and Lemma 1.31]{Di-Pietro.Droniou:20} for a generalization to polytopal elements) to pass to the second line.

  Take now $v\in W^{q+1,s}(T)$ and $0\le m\le q+1$, and apply the discrete inverse inequality \cite[Corollary 1.29]{Di-Pietro.Droniou:20} to $\IRTNT{\ell} (v-\lproj{\ell-1}{T}v)\in\Poly{\ell}(T)^d$ to write
  \begin{align}
    &\seminorm{W^{m,s}(T)^d}{\IRTNT{\ell} (v-\lproj{\ell-1}{T}v)}
    \\
    &\quad
    \begin{aligned}[t]
      &\lesssim h_T^{-m}\norm{L^s(T)^d}{\IRTNT{\ell} (v-\lproj{\ell-1}{T}v)}\nonumber\\
      \overset{\eqref{eq:IRTNT:boundedness}}&\lesssim
      h_T^{-m}\norm{L^s(T)^d}{v-\lproj{\ell-1}{T}v}+h_T^{1-m}\norm{L^s(T)^{d\times d}}{\nabla (v-\lproj{\ell-1}{T}v)}\nonumber\\
      &\lesssim
      h_T^{q+1-m}\seminorm{W^{q+1,s}(T)^d}{v},
      \label{eq:bound.Wmp.IRTNT}
    \end{aligned}
  \end{align}
  where the conclusion follows from the approximation properties of $\lproj{\ell-1}{T}$ in the $L^s$-norm and $W^{1,s}$-seminorm, see \cite[Theorem 1.45]{Di-Pietro.Droniou:20}.
  The approximation property \eqref{eq:IRTNT:approximation} then easily follows by inserting $\pm\lproj{\ell-1}{T}v$ and using $\IRTNT{\ell}\lproj{\ell-1}{T}v=\lproj{\ell-1}{T}v$ (since $\lproj{\ell-1}{T}v\in\RTN{\ell}(T)$) into the seminorm to write
  $$
  \seminorm{W^{m,s}(T)^d}{v - \IRTNT{\ell} v}
  \le \seminorm{W^{m,s}(T)^d}{v - \lproj{\ell-1}{T}v}
  +\seminorm{W^{m,s}(T)^d}{\IRTNT{\ell} (\lproj{\ell-1}{T}v - v)},
  $$
  and conclude thanks to the approximation properties of $\lproj{\ell-1}{T}$ for the first addend and \eqref{eq:bound.Wmp.IRTNT} for the second.

  To prove the trace approximation property \eqref{eq:IRTNT:approximation.trace}, we follow the same argument as in \cite[Theorem 1.45]{Di-Pietro.Droniou:20} by using a continuous trace inequality \cite[Lemma 1.31]{Di-Pietro.Droniou:20} to write
  \begin{multline*}
    \norm{L^s(F)^d}{\partial^\alpha(v - \IRTNT{\ell} v)}
    \\
    \lesssim
    h_T^{-\frac1s}\norm{L^s(T)^d}{\partial^\alpha(v - \IRTNT{\ell} v)}
    +  h_T^{1-\frac1s}  \seminorm{W^{1,s}(T)^d}{\partial^\alpha(v - \IRTNT{\ell} v)},
  \end{multline*}
  and conclude using \eqref{eq:IRTNT:approximation} with $m=r$ and $m=r+1$.
\end{proof}


\section{A scheme with hybrid velocity and pressure}\label{sec:scheme}

In this section we present an extension to the unsteady Navier--Stokes equations of the method for the steady Stokes problem originally introduced in \cite[Section~5.4]{Botti.Botti.ea:24}. This new method hinges on a trilinear form inspired by \cite{Di-Pietro.Krell:18} (see also~\cite{Botti.Di-Pietro.ea:19*1} for a variant) and convective stabilization term in the spirit of \cite{Brezzi.Marini.ea:04}.

\subsection{Discrete spaces and interpolators}\label{sec:discrete.spaces}

Let an integer $k \ge 0$ be fixed and set
\begin{multline*}
  \underline{U}_h^k
  \coloneq \Big\{
  \underline{v}_h = ( (v_T)_{T \in \Th}, (v_F)_{F \in \Fh} )
  \st
  \\
  \text{$v_T \in \RTN{k+1}(T)$ for all $T \in \Th$
    and
    $v_F \in \Poly{k}(F)^d$ for all $F \in \Fh$}
  \Big\}
\end{multline*}
and
\begin{multline*}
  \underline{P}_h^k
  \coloneq \Big\{
  \underline{q}_h = ( (q_T)_{T \in \Th}, (q_F)_{F \in \Fh} )
  \st
  \\
  \text{$q_T \in \Poly{k}(T)$ for all $T \in \Th$
    and
    $q_F \in \Poly{k}(F)$ for all $F \in \Fh$}
  \Big\}.
\end{multline*}
The meaning of the components in the discrete spaces is provided by the interpolators $\IUh : H^1(\Omega)^d \to \Uh$
and $\IPh : H^1(\Omega) \to \Ph$ such that, for all $v \in H^1(\Omega)^d$ and all $q \in H^1(\Omega)$,
\begin{equation}\label{eq:IUh}
  \IUh v
  \coloneq (
  (\IRTNT{k+1} v)_{T \in \Th}, (\lproj{k}{F} v)_{F \in \Fh}
  ),\qquad
  \IPh q
  \coloneq (
  (\lproj{k}{T} q)_{T \in \Th}, (\lproj{k}{F} q)_{F \in \Fh}
  ).
\end{equation}
The main difference with respect to the classical Hybrid High-Order interpolator lies in the fact that $\IRTNT{k+1}$ replaces the $L^2$-orthogonal projector on the element space.
The restrictions to a mesh element $T \in \Th$ of the discrete spaces, of their elements, and of the interpolators are denoted replacing the subscript $h$ with $T$.

Given $\underline{v}_h \in \Uh$ and $\underline{q}_h \in \Ph$, we denote by $v_h : \Omega \to \Real^d$ and $q_h : \Omega \to \Real$ the piecewise polynomial functions such that
\begin{equation}\label{eq:piecewise.functions}
  \text{%
    $(v_h)_{|T} \coloneq v_T$ and $(q_h)_{|T} \coloneq q_T$ for all $T \in \Th$.
  }
\end{equation}
The discrete velocity and pressure are sought in the following spaces, respectively incorporating the boundary condition \eqref{eq:strong:bc} on the velocity and the zero-average condition \eqref{eq:strong:p.zero-average} on the pressure:
\[
\UhZ \coloneq \bigg\{
\underline{v}_h \in \Uh \st \text{$v_F = 0$ for all $F \in \Fhb$}
\bigg\},
\qquad
\PhZ \coloneq \bigg\{
\underline{q}_h \in \Ph \st \int_\Omega q_h = 0
\bigg\}.
\]

\subsection{Viscous term}

Let $T \in \Th$.
We define the velocity reconstruction $\RT : \UT \to \Poly{k+1}(T)^d$ such that, for all $\underline{v}_T \in \UT$,
\begin{equation}\label{eq:pT}
  \begin{gathered}
    \int_T \nabla \RT \underline{v}_T : \nabla w
    = - \int_T v_T \cdot \Delta w
    + \sum_{F \in \FT} \int_F v_F \cdot (\nabla w \, n_{TF})
    \qquad \forall w \in \Poly{k+1}(T)^d,
    \\
    \int_T \RT \underline{v}_T \coloneq
    \begin{cases}
      \sum_{F \in \FT} \frac{d_{TF}}{d} \int_F v_F & \text{if $k = 0$},
      \\
      \int_T v_T & \text{otherwise},
    \end{cases}
  \end{gathered}
\end{equation}
where, for any $F \in \FT$, $d_{TF}$ denotes the distance of the center of mass $\overline{x}_T \coloneq \frac{1}{|T|} \int_T x$ of $T$ from the plane containing $F$.

\begin{remark}[Link with the modified elliptic projector]
  Let $v \in H^1(T)^d$.
  Writing \eqref{eq:pT} for $\underline{v}_T = \IUT v$,
  noticing that $\lproj{k-1}{T} \circ \IRTNT{k+1} = \lproj{k-1}{T}$ by \eqref{eq:IRTNT} with $\ell = k + 1$,
  removing the projectors from the right-hand sides using their definitions,
  and integrating by parts, we obtain
  \[
  \begin{gathered}
    \int_T \nabla \RT \IUT v : \nabla w = \int_T \nabla v : \nabla w
    \qquad \forall w \in \Poly{k+1}(T)^d,
    \\
    \int_T \RT \IUT v = \begin{cases}
      \sum_{F \in \FT} \frac{d_{TF}}{d} \int_F v & \text{if $k = 0$},
      \\
      \int_T v & \text{otherwise}.
    \end{cases}
  \end{gathered}
  \]
  Observe next that, for all $\underline{v}_T \in \underline{U}_T^0$,
  \begin{multline*}
    \int_T \RTo \underline{v}_T
    = \frac1d \int_T \RTo \underline{v}_T \operatorname{div}(x - \overline{x}_T)
    \\
    =
    -\frac1d \cancel{%
      \int_T \nabla \RTo \underline{v}_T \cdot (x - \overline{x}_T)
    }
    + \frac1d \sum_{F \in \FT} \int_F \RTo \underline{v}_T \, (x - \overline{x}_T) \cdot n_{TF}
    = \sum_{F \in \FT} \frac{d_{TF}}{d} \int_F \RTo \underline{v}_T,
  \end{multline*}
  where we have used the fact that $\operatorname{div}(x - \overline{x}_T) = d$ in the first equality,
  an integration by parts together with the facts that $\nabla \RTo \underline{I}_{U,T}^0 v$ is constant inside $T$ and that the vector-valued function $T \ni x \mapsto (x - \overline{x}_T) \in \Real^d$ is $L^2(T)^d$-orthogonal to constant fields by definition of $\overline{x}_T$ in the cancellation,
  and concluded observing that, for all $F \in \FT$, $d_{TF} = (x - \overline{x}_T) \cdot n_{TF}$ for all $x \in F$.
  Accounting for the previous relation, the closure condition for $k = 0$ can be rewritten
  \[
  \sum_{F \in \FT} \frac{d_{TF}}{d} \int_F (\RTo \underline{I}_{U,T}^0 v - v) = 0,
  \]
  showing  that $\RT \circ \IUT$ is in fact the modified elliptic projector of \cite[Section~5.1.2]{Di-Pietro.Droniou:20} with weights $\omega_{TF} = \frac{d_{TF}}{d}$ for all $F \in \FT$.
  As a result, by \cite[Theorem~5.7]{Di-Pietro.Droniou:20} we have the following approximation properties:
  For all $T \in \Th$ and all $v \in H^{k+2}(T)^d$,
  \[ 
  \norm{L^2(\partial T)^{d\times d}}{\nabla \RT \IUT v - \nabla v}
  \lesssim h_T^{k+\frac12} \seminorm{H^{k+2}(T)^d}{v}.
  \] 
  The above operator also corresponds to the Virtual Element projection $\Pi^\nabla_{T,k+1}$ of \cite{Beirao-da-Veiga.Brezzi.ea:14} with a particular choice of the ${\cal R}$ function in \cite[Section 2.2]{Beirao-da-Veiga.Lovadina.ea:17*1}.
\end{remark}

The viscous bilinear form $a_h : \Uh \times \Uh \to \Real$ is such that, for all $(\underline{w}_h, \underline{v}_h) \in \Uh \times \Uh$,
\[
\begin{gathered}
  a_h(\underline{w}_h, \underline{v}_h)
  \coloneq \sum_{T \in \Th} a_T(\underline{w}_T, \underline{v}_T),
  \\
  a_T(\underline{w}_T, \underline{v}_T)
  \coloneq \int_T \nabla \RT \underline{w}_T : \nabla \RT \underline{v}_T
  + s_T(\underline{w}_T, \underline{v}_T).
\end{gathered}
\]
Above, $s_T : \UT \times \UT \to \Real$ denotes a stabilization bilinear form that penalizes the components of $(\delta_T^k \underline{v}_T, (\delta_{TF}^k \underline{v}_T)_{F \in \FT} ) \coloneq \IUT \RT \underline{v}_T - \underline{v}_T$.
Examples of such bilinear form are
\[
s_T(\underline{w}_T, \underline{v}_T)
= \lambda_T h_T^{-2} \int_T \delta_T^k \underline{w}_T \cdot \delta_T^k \underline{v}_T
+ h_T^{-1} \sum_{F \in \FT} \int_F \delta_{TF}^k \underline{w}_T \cdot \delta_{TF}^k \underline{v}_T
\]
with, e.g., $\lambda_T \coloneq \operatorname{card}(\FT) \frac{h_T^d}{|T|}$ to equilibrate the two contributions, or
\[
s_T(\underline{w}_T, \underline{v}_T)
= h_T^{-1} \sum_{F \in \FT} \int_F (\delta_T^k - \delta_{TF}^k) \underline{w}_T \cdot (\delta_T^k - \delta_{TF}^k) \underline{v}_T.
\]
Both choices lead, for all $T \in \Th$, to the uniform seminorm equivalence:
\begin{equation}\label{eq:sT:seminorm.equivalence}
  a_T(\underline{v}_T, \underline{v}_T)^{\frac12} \simeq \norm{1,T}{\underline{v}_T}
  \qquad \forall \underline{v}_T \in \UT,
\end{equation}
with local discrete $H^1$-seminorm
\begin{equation}\label{eq:norm.1.T}
  \norm{1,T}{\underline{v}_T}
  \coloneq \left(\norm{L^2(T)^{d \times d}}{\nabla v_T}^2
  + h_T^{-1} \sum_{F \in \Fh} \norm{L^2(F)^d}{v_F - v_T}^2
  \right)^{\frac12}.
\end{equation}
Moreover, it can be proved using standard arguments (see, e.g., \cite[Proposition~2.14]{Di-Pietro.Droniou:20}) along with the approximation properties of $\IRTNT{k+1}$ that it holds, for all $T \in \Th$,
\begin{equation}\label{eq:sT:consistency}
  s_T(\IUT v, \IUT v)^{\frac12} \lesssim h_T^{k+1} \seminorm{H^{k+2}(T)^d}{v}
  \qquad \forall v \in H^{k+2}(T)^d.
\end{equation}

\begin{remark}[Local stabilization]
  The first choice above for the form $s_T(\cdot,\cdot)$ is in the spirit of the so called ``dofi-dofi'' stabilization of Virtual Elements \cite{Beirao-da-Veiga.Brezzi.ea:13,Beirao-da-Veiga.Lovadina.ea:17*1} while the second one is in the spirit of Hybrid High-Order methods \cite{Di-Pietro.Ern.ea:14,Di-Pietro.Ern:15}.
\end{remark}

\subsection{Velocity--pressure coupling}

Define the pressure gradient $\GT : \PT \to \RTN{k+1}(T)$ such that, for all $\underline{q}_T \in \PT$,
\begin{equation}\label{eq:GT}
  \int_T \GT \underline{q}_T \cdot w
  = - \int_T q_T (\nabla \cdot w)
  + \sum_{F \in \FT} \int_F q_F (w \cdot n_{TF})
  \qquad\forall w \in \RTN{k+1}(T).
\end{equation}
This gradient reconstruction is defined so that the following polynomial consistency property holds:
\begin{equation}\label{eq:GT:polynomial.consistency}
  \GT \IPT q = \nabla q \qquad \forall q \in \Poly{k+1}(T).
\end{equation}
The velocity--pressure coupling is based on the bilinear form $b_h : \Uh \times \Ph \to \Real$ such that, for all $(\underline{v}_h, \underline{q}_h) \in \Uh \times \Ph$,
\begin{equation}\label{eq:bh}
  \text{$b_h(\underline{v}_h, \underline{q}_h)
    \coloneq \sum_{T \in \Th} \int_T v_T \cdot \GT \underline{q}_T$}.
\end{equation}

In what follows, we denote the discrete subspace of divergence-free functions by
\begin{equation}\label{eq:ZhZ}
  \ZhZ \coloneq
  \Big\{ \underline{v}_h \in \UhZ \st
  \text{$b_h(\underline{v}_h, \underline{q}_h)=0$ for all $\underline{q}_h \in \Ph$}
  \Big\}.
\end{equation}

\begin{proposition}[Characterization of $\ZhZ$]\label{prop:ZhZ:characterization}
  It holds
  \begin{multline}\label{eq:ZhZ:characterization}
    \ZhZ =
    \Big\{
    \underline{v}_h \in \UhZ \st
    \text{$\nabla \cdot v_T = 0$ for all $T \in \Th$,
      $v_T \cdot n = 0$ for all $F \in \Fhb$},
    \\
    \text{and $v_{T_1} \cdot n_{T_1F} + v_{T_2} \cdot n_{T_2F} = 0$ for all $F \in \Fh \setminus \Fhb$ shared by $T_1 \neq T_2 \in \Th$}
    \Big\}.
  \end{multline}
  As a consequence, for all $\underline{v}_h \in \ZhZ$, the function $v_h$ defined by \eqref{eq:piecewise.functions} belongs to $\Hdiv{\Omega}$, with $\nabla \cdot v_h = 0$ in $\Omega$ and $v_h \cdot n = 0$ on $\partial\Omega$.
\end{proposition}

\begin{remark}[Element components in $\ZhZ$]
  An important observation is the following:
  \begin{equation}\label{eq:ZhZ.degree.T}
    \text{
      For all $T \in \Th$ and $\underline{v}_h \in \ZhZ$,
      it holds $v_T \in \Poly{k}(T)^d$,
    }
  \end{equation}
  which follows from $v_T \in \RTN{k+1}(T)$ and $\nabla \cdot v_T = 0$.
\end{remark}

\begin{proof}[Proof of Proposition \ref{prop:ZhZ:characterization}]
  Let $\underline{v}_h \in \ZhZ$.
  Expanding first $b_h$ and then $\GT$ according to the respective definitions \eqref{eq:bh} and \eqref{eq:GT}, the condition $b_h(\underline{v}_h, \underline{q}_h) = 0$ translates to
  \begin{equation}\label{eq:discrete:mass:bis}
    \sum_{T \in \Th}
    \left(
    \int_T (\nabla \cdot v_T) q_T
    - \sum_{F \in \FT} \int_F (v_T \cdot n_{TF}) q_F
    \right)
    = 0
    \qquad \forall \underline{q}_h \in \Ph.
  \end{equation}
  Let now $T \in \Th$ be fixed and select $\underline{q}_h$ above such that $q_{T'} = 0$ for all $T' \in \Th \setminus \{ T \}$
  and $q_F = 0$ for all $F \in \Fh$.
  Letting $q_T$ span $\Poly{k}(T)$ and observing that $\nabla \cdot v_T \in \Poly{k}(T)$ since $v_T \in \RTN{k+1}(T)$, we get, since $T$ is generic,
  \begin{equation}\label{eq:div.uT=0}
    \nabla \cdot v_T = 0 \qquad \forall T \in \Th.
  \end{equation}
  Let now $F \in \Fh \setminus \Fhb$ be shared by the mesh elements $T_1 \neq T_2 \in \Th$, and take $\underline{q}_h$ in \eqref{eq:discrete:mass:bis} such that
  $q_T = 0$ for all $T \in \Th$
  and $q_{F'} = 0$ for all $F' \in \Fh \setminus \{ F \}$.
  Letting $q_F$ span $\Poly{k}(F)$ and observing that $v_{T_i} \cdot n_{T_iF} \in \Poly{k}(F)$ for $i \in \{1, 2\}$ (see \cite[Lemma~3.6]{Gatica:14}), we obtain, since $F$ is generic,
  \begin{equation}\label{eq:uT1.nT1F=uT2.nT2F}
    v_{T_1} \cdot n_{T_1F} + v_{T_2} \cdot n_{T_2F} = 0 \qquad \forall F \in \Fh \setminus \Fhb.
  \end{equation}
  A similar reasoning for a boundary face $F \in \Fhb$ gives
  \begin{equation}\label{eq:uT.n=0}
    v_T \cdot n = 0 \qquad \forall F \in \Fhb.
  \end{equation}
  Conversely, since the left-hand side of \eqref{eq:discrete:mass:bis} corresponds to $b_h(\underline{v}_h, \underline{q}_h)$, it can easily be checked that \eqref{eq:div.uT=0}, \eqref{eq:uT1.nT1F=uT2.nT2F}, and \eqref{eq:uT.n=0} imply $b_h(\underline{u}_h, \underline{q}_h) =0$ for all $\underline{q}_h \in \Ph$.
\end{proof}

\subsection{Convective term}\label{sec:scheme:convective.term}

The convective trilinear form $t_h : [ \Uh ]^3 \to \Real$ is such that, for all $(\underline{w}_h, \underline{v}_h, \underline{z}_h) \in [ \Uh ]^3$,
\begin{subequations}\label{eq:th}
  \begin{gather}
    t_h(\underline{w}_h, \underline{v}_h, \underline{z}_h)
    \coloneq \sum_{T \in \Th} t_T(\underline{w}_T, \underline{v}_T, \underline{z}_T),
    \\ \label{eq:tT}
    t_T(\underline{w}_T, \underline{v}_T, \underline{z}_T)
    \coloneq
    \int_T (w_T \cdot \nabla) v_T \cdot z_T
    + \frac12 \sum_{F \in \FT} \int_F (w_T \cdot n_{TF}) (v_F - v_T) \cdot (z_F + z_T).
  \end{gather}
\end{subequations}
This trilinear form can be interpreted as a centered-flux approximation of the convective term and, as such, it satisfies the non-dissipativity property stated in the following lemma.

\begin{lemma}[Non-dissipativity of the convective trilinear form]
  \label{lem:non-diss}
  For all $\underline{w}_h \in \ZhZ$, it holds
  \begin{equation}\label{eq:th:non-dissipativity}
    t_h(\underline{w}_h, \underline{v}_h, \underline{v}_h) = 0
    \qquad \forall \underline{v}_h \in \Uh.
  \end{equation}
\end{lemma}

\begin{proof}
  We write, for all $T \in \Th$,
  \begin{equation}\label{eq:th:non-dissipativity:1}
    \begin{aligned}
      &t_T(\underline{w}_T, \underline{v}_T, \underline{v}_T)
      \\
      &\quad
      \begin{aligned}[t]
        \overset{\eqref{eq:tT}}&=
        \frac 12 \int_T (w_T \cdot \nabla) v_T \cdot  v_T
        + \frac12 \int_T (w_T \cdot \nabla) v_T \cdot  v_T
        \nonumber\\
        &\quad
        + \frac12 \sum_{F \in \FT} \int_F (w_T \cdot n_{TF}) (v_F - v_T) \cdot  (v_F + v_T)
        \nonumber\\
        &= - \frac12 \int_T  v_T \cdot \left[
          \nabla \cdot (  v_T \otimes w_T )
          \right]
        + \frac12 \sum_{F \in \FT} \int_F (w_T \cdot n_{TF}) ( v_T \cdot  v_T)
        \nonumber\\
        &\quad
        + \frac12 \int_T (w_T \cdot \nabla) v_T \cdot  v_T
        + \frac12 \sum_{F \in \FT} \int_F (w_T \cdot n_{TF}) (v_F - v_T) \cdot  (v_F + v_T)
        \nonumber\\
        &=
        \frac12 \sum_{F \in \FT} \int_F (w_T \cdot n_{TF}) (v_F \cdot  v_F),
      \end{aligned}
    \end{aligned}
  \end{equation}
  where $\otimes$ denotes the tensor product of two vectors; above, we have used an integration by parts on the first term in the second equality
  and the product rule together with the fact that $\nabla \cdot w_T = 0$ (since $\underline{w}_h \in \ZhZ$, cf.~\eqref{eq:ZhZ:characterization}) to write $\int_T  v_T \cdot \left[ \nabla \cdot (  v_T \otimes w_T ) \right] = \int_T (w_T \cdot \nabla) v_T \cdot  v_T$ in the third inequality.
  We now observe that, since $\underline{w}_h \in \ZhZ$, the normal traces $(w_T \cdot n_{TF})_{T\in\Th,\,F\in\FT}$ are continuous across interfaces (see again \eqref{eq:ZhZ:characterization}).
  Thus, summing \eqref{eq:th:non-dissipativity:1} over $T\in\Th$ and using the above remark together with the fact that $v_F$ is single-valued on all faces, we obtain
  \[
  \begin{aligned}
    &t_h(\underline{w}_h, \underline{v}_h, \underline{v}_h)
    \\
    &\quad
    \begin{aligned}[t]
      \overset{\eqref{eq:th}}&= \sum_{T \in \Th} t_T(\underline{w}_T, \underline{v}_T, \underline{v}_T)
      \overset{\eqref{eq:th:non-dissipativity:1}}= \frac12 \sum_{T \in \Th} \sum_{F \in \FT} \int_F (w_T \cdot n_{TF}) (v_F \cdot  v_F)\\
      &=
      \frac12 \sum_{F\in\Fh\setminus\Fhb}  \int_F \underbrace{ (w_{T_1} \cdot n_{T_1F} + w_{T_2} \cdot n_{T_2F}) }_{=0}
      (v_F \cdot  v_F)
      \\
      &\quad
      +\frac12 \sum_{F\in\Fh\setminus\Fhb}  \int_F \underbrace{(w_{T} \cdot n_{TF})}_{=0} (v_F\cdot v_F)=0.
    \end{aligned}
  \end{aligned}
  \]
\end{proof}

Given a family of strictly positive real numbers $\beta \coloneq ( \beta_T )_{T \in \Th}$, possibly varying in time,
we also define the convective stabilization bilinear form $j_{\beta,h} : \Uh \times \Uh \to \Real$ such that, for all $(\underline{w}_h, \underline{v}_h) \in \Uh \times \Uh$,
\[
\begin{gathered}
  j_{\beta,h}(\underline{w}_h, \underline{v}_h)
  \coloneq \sum_{T \in \Th} j_{\beta,T}(\underline{w}_T, \underline{v}_T),
  \\
  j_{\beta,T}(\underline{w}_T, \underline{v}_T)
  \coloneq \beta_T \sum_{F \in \FT} \int_F (w_F -  w_T) \cdot (v_F -  v_T).
\end{gathered}
\]
This stabilization contains some upwinding, and its role in the derivation of Reynolds-robust estimates will become clear in Lemma~\ref{lem:Err.conv.h:estimate}; see, in particular, the estimate of the term $\term_2(T)$ on a convection-dominated element $T \in \Th$.

\subsection[Discrete velocity L2-product and norm]{Discrete velocity $L^2$-product and norm}\label{sec:scheme:L2.product.norm}

The discretization of the velocity time derivative hinges on the discrete $L^2$-product $(\cdot, \cdot)_{0,h} : \Uh \times \Uh \to \Real$ such that, for all $(\underline{w}_h, \underline{v}_h) \in \Uh \times \Uh$,
\begin{equation}\label{eq:0.h.product}
  \begin{gathered}
    \text{%
      $(\underline{w}_h, \underline{v}_h)_{0,h} \coloneq
      \sum_{T \in \Th} (\underline{w}_T, \underline{v}_T)_{0,T}$,
    }
    \\
    (\underline{w}_T, \underline{v}_T)_{0,T}
    \coloneq \int_T w_T \cdot v_T
    + h_T \sum_{F \in \FT} \int_F (w_F - w_T) \cdot (v_F - v_T).
  \end{gathered}
\end{equation}
The corresponding global and local norms are respectively given by:
For all $\underline{v}_h \in \Uh$,
\begin{equation}\label{eq:norm.0.T.h}
  \text{%
    $\norm{0,h}{\underline{v}_h} \coloneq (\underline{v}_h, \underline{v}_h)_{0,h}^{\frac12}$
    and
    $\norm{0,T}{\underline{v}_T} \coloneq (\underline{v}_T, \underline{v}_T)_{0,T}^{\frac12}$
    for all $T \in \Th$.
  }
\end{equation}

\subsection{Discrete problem}

In what follows, to simplify the notation, we omit the dependence on the time $t$ whenever it can be inferred from the context.
The discrete problem consists in finding a differentiable function $\underline{u}_h : [0, \tF] \to \UhZ$ and a function $\underline{p}_h : (0, \tF\rbrack \to \PhZ$ such that
\begin{equation}\label{eq:discrete:ic}
  \underline{u}_h(0) = \IUh u_0
\end{equation}
and, for all $t \in \lbrack 0,\tF\rbrack$ and all $(\underline{v}_h,\underline{q}_h) \in \UhZ \times \PhZ$,
\begin{subequations}\label{eq:discrete}
  \begin{alignat}{2}\label{eq:discrete:momentum}
    \left(\frac{d \underline{u}_h}{dt}, \underline{v}_h \right)_{0,h}
    + \nu a_h(\underline{u}_h, \underline{v}_h)
    + t_h(\underline{u}_h, \underline{u}_h, \underline{v}_h)
    + j_{\beta,h}(\underline{u}_h, \underline{v}_h)
    + b_h(\underline{v}_h, \underline{p}_h)
    &= \int_\Omega f \cdot v_h,
    \\ \label{eq:discrete:mass}
    -b_h(\underline{u}_h, \underline{q}_h) &= 0.
  \end{alignat}
\end{subequations}

\begin{remark}[Pointwise divergence-free velocity field]\label{rem:div.uh=0}
  We start by noticing that the mass balance \eqref{eq:discrete:mass} actually holds for any $\underline{q}_h \in \Ph$, as can be checked observing that $b_h(\underline{u}_h, \IPh 1) = 0$ since the gradient reconstruction \eqref{eq:GT} vanishes for interpolates of constant fields by \eqref{eq:GT:polynomial.consistency}.
  As a consequence, $\underline{u}_h \in \ZhZ$ and, recalling Proposition~\ref{prop:ZhZ:characterization}, we have that $u_h \in \Hdiv{\Omega}$ and $u_h \cdot n = 0$ on $\partial \Omega$.
\end{remark}

\begin{remark}[Variations with velocities in the Brezzi--Douglas--Marini space]
  Let us investigate some scenarios in which the Brezzi--Douglas--Marini space is used instead of the Raviart-Thomas--Nédélec space for the velocity field.

  \medskip
  \noindent
      {\textsc{Choice 1.}}
      The first case is obtained considering the following discrete spaces:
      \[
      \begin{gathered}
        \Uh \coloneqq \left( \bigtimes_{T \in \Th} \mathcal{BDM}^{k+1}(T) \right) \times \left( \bigtimes_{F \in \Fh} \Poly{k}(F)^d \right),
        \\
        \Ph \coloneqq \left( \bigtimes_{T \in \Th} \Poly{k}(T) \right) \times \left( \bigtimes_{F \in \Fh} \Poly{k+1}(F) \right).
      \end{gathered}
      \]
      We generalize the definitions of the trilinear form \eqref{eq:tT} and the convective stabilization bilinear form as follows:
      \[
      \begin{aligned}
        t_T (\underline{w}_T, \underline{v}_T, \underline{z}_T)
        &\coloneqq
        \int_T (w_T \cdot \nabla) \lproj{k}{T}v_T \cdot \lproj{k}{T}z_T
        \\
        &\quad
        + \frac12 \sum_{F \in \FT} \int_F (w_T \cdot n_{TF}) (v_F - \lproj{k}{T}v_T) \cdot (z_F + \lproj{k}{T}z_T), \\
        j_{\beta,T}(\underline{w}_T, \underline{v}_T)
        &\coloneqq
        \beta_T \sum_{F \in \FT} \int_F (w_F -  \lproj{k}{T}w_T) \cdot (v_F -  \lproj{k}{T}v_T).
      \end{aligned}
      \]
      Exploiting the approximation properties of the interpolator ${I}_{\mathcal{BDM},T}^{k+1}$, and under additional regularity on $(u \cdot \nabla) u$, one can prove that this method attains the same convergence rate stated in Corollary \ref{cor:convergence.rate} below.
      Modifying the convective forms by adding the projectors $\lproj{k}{T}$ is crucial to avoid losing one order of convergence for the convective error component. This loss would occur because the key property \eqref{eq:ZhZ.degree.T}, used in Lemma \ref{lem:Err.conv.h:estimate}, does not hold for BDM functions.

      \medskip
      \noindent
          {\textsc{Choice 2.}}
          In the second case, we keep the same discrete operators and define:
          \[
          \begin{gathered}
            \Uh \coloneqq \left( \bigtimes_{T \in \Th} \mathcal{BDM}^{k}(T) \right) \times \left( \bigtimes_{F \in \Fh} \Poly{k}(F)^d \right),
            \\
            \Ph \coloneqq \left( \bigtimes_{T \in \Th} \Poly{k-1}(T) \right) \times \left( \bigtimes_{F \in \Fh} \Poly{k}(F) \right).
          \end{gathered}
          \]
          This choice of function spaces coincides with that in  \cite[Section~5.3]{Botti.Botti.ea:24} for the Stokes problem, where the consistency error of the diffusion term is of order $O(h^k)$. However, the viscous bilinear form used here differs from the one adopted there. Moreover, using
          the approximation property of $I_{\mathcal{BDM},T}^{k}$, the identity $\lproj{k-1}{T} I_{\mathcal{BDM},T}^{k}v = \lproj{k-1}{T} v$ for divergence-free $v \in H^1(T)^d$, and standard arguments \cite[Lemma 2.18]{Di-Pietro.Droniou:20},
          one can achieve the same convergence rates up to $h^{k+1}$ as in the method considered throughout this paper.
          Actually, it can also be proved that the discrete velocity resulting from this method corresponds to that of the original scheme, i.e., the scheme \eqref{eq:discrete:ic} with spaces as defined in Section \ref{sec:discrete.spaces}.
          We finally observe that such faster convergence rate in diffusion dominated cases compares favorably also with the method of \cite{Rhebergen.Wells:18} (in this respect, see also Remark \ref{rem:keegan.et.al}).
          Notice that a higher convergence rate in the diffusion-dominated regime compared to the convection-dominated regime means that the accuracy of the method improves instead of getting worse as we refine the mesh.
\end{remark}

\subsection{Well-posedness of the discrete problem}

The main theorem of this section states the existence and uniqueness of the solution to the scheme, together with a priori estimates on the discrete velocity. These estimates are established in the norm $\tnorm{h}{{\cdot}}$ on $C^0([0,\tF];\UhZ)$ defined by: For all $\underline{v}_h\in C^0([0,\tF];\UhZ)$,
\begin{equation}\label{eq:def:tnorm}
  \tnorm{h}{\underline{v}_h}^2 \coloneq
  \max_{t \in [0,\tF]} \norm{0,h}{\underline{v}_h(t)}^2
  + \int_0^{\tF} \nu \norm{1,h}{\underline{v}_h(\tau)}^2 \, d \tau
  + \int_0^{\tF} \seminorm{\beta,h}{\underline{v}_h(\tau)}^2 \, d \tau,
\end{equation}
with diffusive $H^1$-like seminorm
\begin{equation}\label{eq:norm.1.h}
  \text{
    $\norm{1,h}{\underline{v}_h} \coloneq \left(
    \sum_{T \in \Th} \norm{1,T}{\underline{v}_T}^2
    \right)^{\frac12}$
    with $\norm{1,T}{{\cdot}}$ given by \eqref{eq:norm.1.T}
  }
\end{equation}
and convective seminorm $\seminorm{\beta,h}{{\cdot}}$ such that, for all $\underline{v}_h \in \Uh$,
\begin{equation}\label{eq:seminorm.beta}
  \text{
    $\seminorm{\beta,h}{\underline{v}_h} \coloneq j_{\beta,h}(\underline{v}_h, \underline{v}_h)^{\frac12}$
    and $\seminorm{\beta,T}{\underline{v}_T} \coloneq j_{\beta,T}(\underline{v}_T, \underline{v}_T)^{\frac12}$ for all $T \in \Th$.
  }
\end{equation}

\begin{remark}[Choice of norm]
  The choice \eqref{eq:def:tnorm} is tailored to the analysis of the scheme and inspired by the norm that naturally appears in the well-posedness study of the continuous problem.
  Specifically, the first term in \eqref{eq:def:tnorm} plays, at the discrete level, the role of the $C^0([0,\tF]; L^2(\Omega)^d)$-norm,
  the second term is akin to an $L^2(0,\tF; H_0^1(\Omega)^d)$-norm,
  while the last one results from the addition of convective stabilization
  and, although not having a direct continuous counterpart, it is associated to control on $|\beta\cdot \nabla v|$ in a DG sense.
\end{remark}

\begin{lemma}[Inf-sup condition on $b_h$]
  \label{lem:inf-sup}
  It holds
  \[
  \left(
  \norm{L^2(\Omega)}{q_h}^2
  + \sum_{T \in \Th} h_T^2 \norm{L^2(T)^d}{\GT \underline{q}_T}^2
  \right)^{\frac12}
  \lesssim
  \sup_{\underline{v}_h \in \UhZ \setminus \{ \underline{0} \}}
  \frac{b_h( \underline{v}_h, \underline{q}_h)}{\norm{1,h}{\underline{v}_h}}
  \qquad \forall \underline{q}_h \in \PhZ.
  \]
\end{lemma}
\begin{proof}
  Straightforward consequence of \cite[Lemma~3 and Theorem~14]{Botti.Botti.ea:24}.
\end{proof}

\begin{theorem}[Well-posedness of the scheme]\label{thm:well.posedness.scheme}
  Assuming that $\beta_T\in C^0([0,\tF])$ for all $T\in\Th$, there is a unique solution $(\underline{u}_h,\underline{p}_h)\in C^1([0,\tF];\UhZ)\times C^0([0,\tF];\PhZ)$ to \eqref{eq:discrete}, which moreover satisfies
  \begin{equation}\label{eq:estimate.uh}
    \tnorm{h}{\underline{u}_h}^2\lesssim e^{\tF}\left(\int_0^{\tF} \norm{L^2(\Omega)^d}{f(\tau)}^2\,d\tau
    + \norm{H^1(\Omega)^d}{u_0}^2\right).
  \end{equation}
\end{theorem}

\begin{proof}
  The scheme \eqref{eq:discrete} shows that $\underline{u}_h:[0,\tF]\to\ZhZ$ solves the ordinary differential equation \eqref{eq:discrete:momentum} in which the test functions $\underline{v}_h$ are taken in $\ZhZ$ (which cancels the pressure term $b_h(\underline{v}_h,\underline{p}_h))$.
  This ODE is a nonlinear equation in a finite dimensional space with continuous coefficients (due to the assumption on $\beta_T$), and therefore has a unique $C^1$ local-in-time solution.
  Assuming that we can prove the estimate \eqref{eq:estimate.uh}, we see that $\underline{u}_h$ does not blow up in finite time, and therefore exists up to $\tF$. The pressure $\underline{p}_h$ is then recovered by writing \eqref{eq:discrete:momentum} for any test function $\underline{v}_h$ which, by Lemma \ref{lem:inf-sup}, uniquely defines this pressure as a linear map of the other terms in the equation; since all these terms are continuous in time, this shows that $\underline{p}_h\in C^0([0,\tF];\PhZ)$.

  It remains to show \eqref{eq:estimate.uh}. Making $\underline{v}_h=\underline{u}_h$ in \eqref{eq:discrete:momentum} and using the non-dissipativity \eqref{eq:th:non-dissipativity} and the norm equivalence \eqref{eq:sT:seminorm.equivalence}, we obtain
  \[
  \frac{d}{dt}\norm{0,h}{\underline{u}_h}^2+\nu\norm{1,h}{\underline{u}_h}^2+\seminorm{\beta,h}{\underline{u}_h}^2\lesssim
  \norm{L^2(\Omega)^d}{f}\norm{L^2(\Omega)^d}{u_h}.
  \]
  Integrating over $[0,t]$ for an arbitrary $t\in [0,\tF]$, using Young's inequality on the right-hand side together with $\norm{L^2(\Omega)^d}{u_h}\le \norm{0,h}{\underline{u}_h}$, and finally invoking the Gronwall inequality \cite[Proposition 2.1]{Emmrich:99} as well as the fact that $\norm{0,h}{\underline{u}_h(0)}=\norm{0,h}{\IUh u_0} \lesssim \norm{H^1(\Omega)^d}{u_0}$ concludes the proof.
\end{proof}

\begin{remark}[Comparison with the method of Rhebergen and Wells]\label{rem:keegan.et.al}
  After discretization in time, linearization, and static condensation of the element unknowns (with the possible exception of one pressure unknown per element), the discrete problem~\eqref{eq:discrete} translates into linear systems with analogous size and sparsity pattern as for the method originally introduced in \cite{Rhebergen.Wells:18} using the Crank--Nicolson method to advance in time and later analyzed in \cite{Keegan.Horvath.ea:23} using more general time-stepping strategies.
  However, the  scheme \eqref{eq:discrete} converges in space as $h^{k+1}$ as opposed to $h^k$ for \cite{Rhebergen.Wells:18,Keegan.Horvath.ea:23} in the diffusion-dominated regime, while still being pressure-robust.
  This results from two important differences: the use of a slightly larger space for the element velocity ($\RTN{k+1}(T)$ as opposed to $\mathcal{BDM}^k(T)$) and of a high-order viscous stabilization, a crucial point being the consistency property \eqref{eq:sT:consistency} of the latter; cf.\ for instance \cite{Cockburn.Di-Pietro.ea:16} for a deep investigation on this subject and \cite{Beirao-da-Veiga.Brezzi.ea:14} for such stabilization in the Virtual Elements framework.
  The present analysis moreover highlights the Reynolds-semi-robustness of the method, and accounts for pre-asymptotic convergence rates through a local Reynolds number.
  Finally, note that the results in \cite{Keegan.Horvath.ea:23} are convection quasi-robust only under a strict data assumption, actually implying $\| u \|_{L^\infty(0,T;H^1(\Omega)^d)} \lesssim \nu$ (therefore not really covering convection dominated cases).
  We are confident that our techniques can be applied to the method of \cite{Keegan.Horvath.ea:23} (which could, in passing, lead to the weakening of certain assumptions) and, conversely, that the time stepping technique proposed therein could be adapted to the present setting.
\end{remark}


\section{Error analysis}\label{sec:error-analysis}

In the present section we develop a convergence analysis for the proposed method, highlighting, in particular, the pressure-robust and Reynolds-semi-robust nature of the scheme.

\subsection{Preliminary results}

In addition to the local and global discrete $L^2$-norms defined by \eqref{eq:norm.0.T.h} and the local and global $H^1$-seminorms respectively defined by \eqref{eq:norm.1.T} and \eqref{eq:norm.1.h}, we will need the following discrete $W^{1,\infty}$-norm for the velocity:
For all $\underline{v}_h \in \UhZ$,
\begin{equation}\label{eq:norm.1.infty}
  \begin{gathered}
    \norm{1,\infty,h}{\underline{v}_h} \coloneq
    \max_{T \in \Th} \norm{1,\infty,T}{\underline{v}_T},
    \\
    \norm{1,\infty,T}{\underline{v}_T}
    \coloneq \norm{L^\infty(T)^{d \times d}}{\nabla v_T}
    + h_T^{-1} \max_{F \in \FT} \norm{L^\infty(F)^d}{v_F - v_T}.
  \end{gathered}
\end{equation}

\begin{lemma}[Boundedness of the convective trilinear form]\label{lem:estimate.th}%
  For all $(\underline{w}_h, \underline{v}_h, \underline{z}_h) \in [ \Uh ]^3$, it holds,
  \begin{equation}\label{eq:estimate.th}
    t_h(\underline{w}_h, \underline{v}_h, \underline{z}_h)
    \lesssim \norm{1,\infty,h}{\underline{v}_h}
    \norm{L^2(\Omega)^d}{w_h}
    \norm{0,h}{\underline{z}_h}.
  \end{equation}
\end{lemma}
\begin{proof}
  Let an element $T \in \Th$ be fixed and let us estimate the terms in the right-hand side of \eqref{eq:tT}.

  Using, respectively, $(2,\infty,2)$- and $(2,\infty,\infty,2)$-H\"{o}lder inequalities in the element and boundary terms that appear in the definition \eqref{eq:tT} of $t_T(\cdot,\cdot,\cdot)$, as well as $\norm{L^\infty(F)^d}{n_{TF}} \le 1$, we get
  \[
  \begin{aligned}
    t_T(\underline{w}_T, \underline{v}_T, \underline{z}_T)
    &\lesssim
    \norm{L^2(T)^d}{w_T}
    \norm{L^\infty(T)^{d\times d}}{\nabla v_T}
    \norm{L^2(T)^d}{ z_T}
    \\
    &\quad
    + \sum_{F \in \FT}
    \norm{L^2(T)^d}{w_T}~
    h_T^{-1}\norm{L^\infty(F)^d}{v_F - v_T}
    \\
    &\qquad
    \times
    \left(
    h_T^{\frac12}\norm{L^2(F)^d}{z_F - z_T}
    + \norm{L^2(T)^d}{z_T}
    \right),
  \end{aligned}
  \]
  where we have additionally used a triangle inequality to write $\norm{L^2(F)^d}{z_F + z_T} \le \norm{L^2(F)^d}{z_F - z_T} + 2 \norm{L^2(F)^d}{z_T}$ together with discrete trace inequalities (see, e.g., \cite[Lemma~1.32]{Di-Pietro.Droniou:20}) to write
  $\norm{L^2(F)^d}{w_T} \lesssim h_T^{-\frac12} \norm{L^2(T)^d}{w_T}$
  and $\norm{L^2(F)^d}{z_T} \lesssim h_T^{-\frac12} \norm{L^2(T)^d}{z_T}$.
  Applying $(2,\infty,2)$-H\"older and Cauchy--Schwarz inequalities on the sums,
  and recalling the definitions \eqref{eq:norm.1.infty} of the $\norm{1,\infty,T}{{\cdot}}$-norm and \eqref{eq:norm.0.T.h} of the $\norm{0,T}{{\cdot}}$-norm, we get
  \[
  t_T(\underline{w}_T, \underline{v}_T, \underline{z}_T)
  \lesssim
  \norm{1,\infty,T}{\underline{v}_T}
  \norm{L^2(T)^d}{w_T}
  \norm{0,T}{\underline{z}_T}.
  \]
  Summing this inequality over $T \in \Th$ and using an $(\infty,2,2)$-H\"{o}lder inequality on the sum, the conclusion follows.
\end{proof}

\begin{lemma}[$W^{1,\infty}$-boundedness of the velocity interpolator]%
  It holds
  \begin{equation}\label{eq:W.1.infty.boundedness.IUh}
    \norm{1,\infty,h}{\IUh w} \lesssim \seminorm{W^{1,\infty}(\Omega)^d}{w}
    \qquad
    \forall w \in W^{1,\infty}(\Omega)^d.
  \end{equation}
\end{lemma}

\begin{proof}
  Let $T \in \Th$. We have
  \begin{align}
    \norm{1,\infty,T}{\IUT w}
    \overset{\eqref{eq:norm.1.infty},\,\eqref{eq:IUh}}&={}
    \seminorm{W^{1,\infty}(T)^d}{\IRTNT{k+1} w}
    + h_T^{-1}\max_{F \in \FT} \norm{L^\infty(F)^d}{\lproj{k}{F} w-\IRTNT{k+1} w }\nonumber\\
    &\lesssim{}  \seminorm{W^{1,\infty}(T)^d}{w}
    + h_T^{-1}\max_{F \in \FT} \norm{L^\infty(F)^d}{\lproj{k}{F} w-\IRTNT{k+1} w },
    \label{eq:W.1.infty.boundedness.IUh:basic}
  \end{align}
  where we have used the boundedness \eqref{eq:boundedness.IRTNT} of $\IRTNT{\ell}$ with $(s,\ell,q)=(\infty,k+1,0)$ to conclude.
  We then write, using triangle inequalities,
  \begin{equation}\label{eq:W.1.infty.boundedness.IUh:T2}
    \begin{aligned}
      &\norm{L^\infty(F)^d}{\lproj{k}{F} w-\IRTNT{k+1} w }
      \\
      &\quad
      \begin{aligned}[t]
        &\le
        \norm{L^\infty(F)^d}{\lproj{k}{F} w-\lproj{k}{T}w}
        +\norm{L^\infty(F)^d}{\lproj{k}{T} w-w}
        +\norm{L^\infty(F)^d}{w-\IRTNT{k+1} w }
        \nonumber\\
        &\lesssim
        \norm{L^\infty(F)^d}{\lproj{k}{T} w-w}
        +\norm{L^\infty(F)^d}{w-\IRTNT{k+1} w }
        \nonumber\\
        &\lesssim
        h_T\seminorm{W^{1,\infty}(T)^d}{w},
      \end{aligned}
    \end{aligned}
  \end{equation}
  where, to pass to the second line, we have written
  \[
  \norm{L^\infty(F)^d}{\lproj{k}{F} w-\lproj{k}{T}w}
  = \norm{L^\infty(F)^d}{\lproj{k}{F} (w-\lproj{k}{T}w)}
  \lesssim \norm{L^\infty(F)^d}{w-\lproj{k}{T}w}
  \]
  by the $L^\infty$-boundedness of the $L^2$-projector \cite[Lemma 1.44]{Di-Pietro.Droniou:20}, and we have concluded using the approximation properties of $\lproj{k}{T}$ (see \cite[Theorem 1.45]{Di-Pietro.Droniou:20}) and of $\IRTNT{k+1}$ (corresponding to \eqref{eq:IRTNT:approximation.trace} with $(s,\ell,q,r)=(\infty,k+1,0,0)$).
  Plugging \eqref{eq:W.1.infty.boundedness.IUh:T2} into \eqref{eq:W.1.infty.boundedness.IUh:basic} concludes the proof.
\end{proof}

\begin{proposition}[Consistency of the velocity--pressure coupling]\label{prop:Interpolate.ZhZ}
  Let $w \in H_0^1(\Omega)^d$ be such that $\nabla \cdot w = 0$.
  Then, $\IUh w \in \ZhZ$, i.e.,
  \begin{equation}\label{eq:bh.IUh.w=0}
    b_h(\IUh w, \underline{q}_h) = 0 \qquad \forall \underline{q}_h \in \PhZ.
  \end{equation}
\end{proposition}

\begin{proof}
  By the definition \eqref{eq:IUh} of the interpolator on $\Uh$
  together with the commutation property \eqref{eq:IRTNT:commutation} for the Raviart--Thomas--Nédélec interpolator,
  $\nabla \cdot \IRTNT{k+1} w = \lproj{k}{T}(\nabla \cdot w)=0$ for all $T \in \Th$.
  Moreover, since $w \in H_0^1(\Omega)^d$, its normal trace is single-valued on interfaces $F \in \Fh \setminus \Fhb$ and it vanishes on boundary faces $F \in \Fhb$.
  Recalling the characterization \eqref{eq:ZhZ:characterization} of $\ZhZ$, this shows that $\IUh w \in \ZhZ$, which is precisely the sought result.
\end{proof}


\subsection{Error estimate and convergence rate for locally smooth solutions}

In this section we state a basic error estimate from which we infer convergence rates for locally smooth solutions.
The estimate of the convective error component provided in Lemma~\ref{lem:Err.conv.h:estimate} below, which is the most subtle, depends on the local regime. For this reason, for any mesh element $T \in \Th$, we define the local Reynolds number
\begin{equation}\label{eq:ReT}
  \ReT \coloneq \frac{
    \big(\beta_T + \norm{L^\infty(T)^d}{u}\big) h_T
  }{\nu}.
\end{equation}
We additionally introduce the time-dependent function
\begin{equation}\label{eq:chi}
  \chi\coloneq\max_{T\in\Th,\,\ReT > 1}\frac{\norm{L^\infty(T)^d}{u}}{\beta_T}.
\end{equation}

\begin{theorem}[Error estimate]\label{thm:error.estimate}
  Assume that $\beta_T\in C^0([0,\tF])$ for all $T\in\Th$.
  Furthermore, let the weak solution to \eqref{eq:strong} be such that $p \in L^2(0,\tF;H^1(\Omega))$ and $u \in H^1(0,\tF; H_0^1(\Omega)^d) \cap L^2(0,\tF; H^2(\Th)^d) \cap L^2(0,\tF;W^{1,\infty}(\Omega)^d)$.
  Set, for the sake of brevity,
  $\huline{u}_h \coloneq \IUh u$,
  $\huline{p}_h \coloneq \IPh p$,
  and define the velocity error
  \begin{equation}\label{eq:eh}
    \underline{e}_h \coloneq \underline{u}_h - \huline{u}_h.
  \end{equation}
  Let the consistency error linear form $\Err_h : \UhZ \to \Real$ be such that, for all $\underline{v}_h \in \ZhZ$,
  \begin{equation}\label{eq:Err.h}
    \Err_h(\underline{v}_h)\coloneq \Err_{{\rm time},h}(\underline{v}_h)+\Err_{{\rm diff},h}(\underline{v}_h)+\Err_{{\rm conv},h}(\underline{v}_h)
  \end{equation}
  with
  \begin{gather}\label{eq:E.time.h}
    \Err_{{\rm time},h}(\underline{v}_h)\coloneq \int_\Omega \frac{d u}{dt} \cdot v_h
    - \left( \frac{d \huline{u}_h}{dt}, \underline{v}_h \right)_{0,h},
    \\ \label{eq:E.diff.h}
    \Err_{{\rm diff},h}(\underline{v}_h)\coloneq  - \int_\Omega \nu \Delta u \cdot v_h - \nu a_h(\huline{u}_h, \underline{v}_h),
    \\ \label{eq:E.conv.h}
    \Err_{{\rm conv},h}(\underline{v}_h)\coloneq \int_\Omega (u \cdot \nabla) u \cdot v_h - t_h(\huline{u}_h, \huline{u}_h, \underline{v}_h)
    - j_{\beta,h}(\huline{u}_h, \underline{v}_h).
  \end{gather}
  Then, recalling the definition \eqref{eq:def:tnorm} of $\tnorm{h}{{\cdot}}$, it holds
  \begin{equation}\label{eq:basic.error.estimate}
    \tnorm{h}{\underline{e}_h}^2  \lesssim K(u) \int_0^{\tF} |\Err_h(\underline{e}_h(\tau))| \, d \tau,
  \end{equation}
  with hidden constant possibly depending on the mesh regularity parameter and the polynomial degree $k$, but independent of the meshsize $h$, of $\tF$, of $\nu$, and of $\beta$, and where
  \[
  K(u) \coloneq \exp\bigg(C\int_0^{\tF}\seminorm{W^{1,\infty}(\Omega)^d}{u(\tau)}\,d\tau\bigg)
  \]
  for a real number $C > 0$ depending only on the mesh regularity parameter and the polynomial degree $k$.
\end{theorem}

\begin{proof}
  We start by noticing that, as a consequence of \eqref{eq:discrete:mass} and of \eqref{eq:bh.IUh.w=0} applied to $w = u$, the vectors of polynomials $\underline{u}_h$, $\huline{u}_h$ and, consequently, their difference $\underline{e}_h$ all belong to the space $\ZhZ$ defined by~\eqref{eq:ZhZ}.
  Moreover, the regularity assumption on $u$ and Theorem \ref{thm:well.posedness.scheme} ensure that $\underline{e}_h \in H^1(0,\tF;\ZhZ)$, which justifies the manipulations below.
  Take $\underline{v}_h = \underline{e}_h$ in the discrete momentum balance equation \eqref{eq:discrete:momentum},
  subtract from the resulting expression the left-hand side with $(\underline{u}_h, \underline{p}_h)$ replaced by $(\huline{u}_h, \huline{p}_h)$,
  and recall that $\underline{u}_h, \huline{u}_h \in \ZhZ$ to cancel the terms involving $\underline{p}_h$, to write, for almost every $t \in (0,\tF)$,
  \[
  \left(\frac{d \underline{e}_h}{dt}, \underline{e}_h\right)_{0,h}
  + \nu a_h(\underline{e}_h, \underline{e}_h)
  + t_h(\underline{u}_h, \underline{u}_h, \underline{e}_h)
  - t_h(\huline{u}_h, \huline{u}_h, \underline{e}_h)
  + j_{\beta,h}(\underline{e}_h, \underline{e}_h)
  = \Err_h(\underline{e}_h),
  \]
  where, in the right-hand side, we have used the fact that $f = \frac{du}{dt} - \nu \Delta u + (u \cdot \nabla) u + \nabla p$ by \eqref{eq:strong:momentum}.
  Notice that no pressure error term appears in the consistency error linear form \eqref{eq:Err.h} because $b_h(\underline{e}_h, \underline{p}_h) = b_h(\underline{e}_h, \widehat{\underline{p}}_h) = 0$ since $\underline{e}_h \in \ZhZ$.
  For this reason, irrotational forcing terms $f$ will lead to vanishing discrete velocity solutions.

  Recalling the definitions \eqref{eq:norm.0.T.h} of the $\norm{0,h}{{\cdot}}$-norm,
  \eqref{eq:seminorm.beta} of the $\seminorm{\beta,h}{{\cdot}}$-seminorm,
  and using the $\norm{1,h}{{\cdot}}$-coercivity of $a_h$ (obtained by summing \eqref{eq:sT:seminorm.equivalence} over $T\in\Th$), we have
  \begin{equation}\label{eq:basic.error.esimate:basic}
    \frac{d \norm{0,h}{\underline{e}_h}^2}{dt}
    + \nu \norm{1,h}{\underline{e}_h}^2
    + \seminorm{\beta,h}{\underline{e}_h}^2
    \lesssim \Err_h(\underline{e}_h)
    + t_h(\huline{u}_h, \huline{u}_h, \underline{e}_h)
    - t_h(\underline{u}_h, \underline{u}_h, \underline{e}_h).
  \end{equation}
  Let us focus on the last two terms in the right-hand side.
  Adding and subtracting $t_h(\underline{u}_h, \huline{u}_h, \underline{e}_h)$ and using the trilinearity of $t_h$, we get
  \[
  \begin{aligned}
    t_h(\huline{u}_h, \huline{u}_h, \underline{e}_h)
    - t_h(\underline{u}_h, \underline{u}_h, \underline{e}_h)
    &= - t_h(\underline{e}_h, \huline{u}_h, \underline{e}_h)
    - \cancel{t_h(\underline{u}_h, \underline{e}_h, \underline{e}_h)},
    \\
    \overset{\eqref{eq:estimate.th}\,\eqref{eq:norm.0.T.h}}&\lesssim
    \norm{1,\infty,h}{\huline{u}_h}
    \norm{0,h}{\underline{e}_h}^2
    \overset{\eqref{eq:W.1.infty.boundedness.IUh}} \lesssim \seminorm{W^{1,\infty}(\Omega)^d}{u}\norm{0,h}{\underline{e}_h}^2,
  \end{aligned}
  \]
  where the cancellation follows from the non-dissipativity \eqref{eq:th:non-dissipativity} of $t_h$.
  Plugging this estimate into \eqref{eq:basic.error.esimate:basic} we get, for any $t \in (0,\tF)$,
  \[
  \frac{d \norm{0,h}{\underline{e}_h}^2}{dt}
  + \nu \norm{1,h}{\underline{e}_h}^2
  + \seminorm{\beta,h}{\underline{e}_h}^2
  \lesssim \Err_h(\underline{e}_h)
  + \seminorm{W^{1,\infty}(\Omega)^d}{u}\norm{0,h}{\underline{e}_h}^2.
  \]
  Integrating in time from $0$ to $t \in (0, \tF)$ and using the fact that $\underline{e}_h(0) = \underline{0}$ by \eqref{eq:eh} together with \eqref{eq:discrete:ic}, we get
  \begin{multline*}
    \norm{0,h}{\underline{e}_h(t)}^2
    + \int_0^t \left(
    \nu \norm{1,h}{\underline{e}_h(\tau)}^2
    + \seminorm{\beta,h}{\underline{e}_h(\tau)}^2
    \right) \, d \tau\\
    \lesssim
    \int_0^t |\Err_h(\underline{e}_h(\tau))| \, d \tau
    + \int_0^t \seminorm{W^{1,\infty}(\Omega)^d}{u(\tau)}
    \norm{0,h}{\underline{e}_h(\tau)}^2 \, d \tau.
  \end{multline*}
  The conclusion then follows from the Gronwall lemma \cite[Proposition 2.1]{Emmrich:99}.
\end{proof}

\begin{corollary}[Convergence rates for smooth solutions]\label{cor:convergence.rate}
  Under the assumptions of Theorem~\ref{thm:error.estimate}, and further assuming that
  $u \in L^2(0,\tF; H^{k+2}(\Th)^d)$ and $\frac{du}{dt} \in L^2(0,\tF; H^{k+1}(\Th)^d)$ for all $T\in\Th$, it holds
  \begin{equation}\label{eq:coroll.rate}
    \begin{aligned}
      \tnorm{h}{\underline{e}_h}
      &\lesssim
      K(u) \Bigg[ \sum_{T \in \Th} h_T^{2(k+1)}
        \int_0^{\tF} \begin{aligned}[t]
          \bigg(
          &\tF \Seminorm{H^{k+1}(T)^d}{\frac{du}{dt}}^2
          + \nu \Seminorm{H^{k+2}(T)^d}{u}^2
          \\
          & + \tF \seminorm{W^{1,\infty}(T)^d}{u}^2 \seminorm{H^{k+1}(T)^d}{u}^2
        \bigg)
        \end{aligned}
        \\
        &\quad
        + \max\left\{ 1, \norm{C^0([0,\tF])}{\chi}\right\}
        \\
        &\qquad
        \times
        \sum_{T \in \Th} h_T^{2k+1} \int_0^{\tF} \min(1,\ReT)
        \big(\beta_T {+} \norm{L^\infty(T)^d}{u} \big)
        \seminorm{H^{k+1}(T)^d}{u}^2
        \Bigg]^{\frac12},
    \end{aligned}
  \end{equation}
  with $K(u)$ as in Theorem~\ref{thm:error.estimate} and hidden constant possibly depending on the mesh regularity parameter and the polynomial degree $k$, but independent of the meshsize $h$, of $\tF$, of $\nu$, and of $\beta$.
\end{corollary}

\begin{remark}[Convergence rate for small Reynolds]\label{rem:convergence.rate}
  If $\ReT \le h_T$ for all $T \in \Th$, then all the terms in the right-hand side of \eqref{eq:coroll.rate} are of order $h^{k+1}$.
\end{remark}

\begin{proof}[Proof of Corollary~\ref{cor:convergence.rate}]
  Note first that, by \eqref{eq:chi}, $u \in C^0([0,\tF]; W^{1,\infty}(\Omega)^d)$ and $\beta_T\in C^0([0,\tF]; (0,\infty))$ for all $T\in\Th$ such that $\ReT>1$ imply $\chi \in C^0([0,\tF])$.
  Plug then the estimates of Lemmas~\ref{lem:estimate.Err.time},~\ref{lem:estimate.Err.diff.h}, and~\ref{lem:Err.conv.h:estimate} below with $\underline{v}_h = \underline{e}_h$ into \eqref{eq:basic.error.estimate} and use Cauchy--Schwarz inequalities to get
  \begin{equation}\label{eq:rate.1}
    \begin{aligned}
      \tnorm{h}{\underline{e}_h}^2 \lesssim{}& K(u) \left(\int_0^{\tF}
      \sum_{T \in \Th} h_T^{2(k+1)}
      \Seminorm{H^{k+1}(T)^d}{\frac{d u}{dt}}^2
      \right)^{\frac12}\left(\int_0^{\tF} \norm{0,h}{\underline{e}_h}^2\right)^{\frac12}
      \\
      &+K(u)\left(\int_0^{\tF}
      \sum_{T \in \Th} \nu h_T^{2(k+1)} \seminorm{H^{k+2}(T)^d}{u}^2
      \right)^{\frac12}\left( \int_0^{\tF} \nu \norm{1,h}{\underline{e}_h}^2\right)^{\frac12}
      \\
      &+K(u) \left(\int_0^{\tF}
      \sum_{T \in \Th} h_T^{2(k+1)}
      \seminorm{W^{1,\infty}(T)^d}{u}^2 \seminorm{H^{k+1}(T)^d}{u}^2
      \right)^{\frac12}\left(\int_0^{\tF} \norm{L^2(\Omega)^d}{e_h}^2\right)^{\frac12}
      \\
      &+K(u) \left( 1 + \norm{C^0([0,\tF]}{\chi}\right)^{\frac12}
      \\
      &\qquad \times
      \left(\int_0^{\tF}
      \sum_{T \in \Th} h_T^{2k+1} \min(1,\ReT)
      \big(
      \beta_T +  \norm{L^\infty(T)^d}{u}
      \big) \seminorm{H^{k+1}(T)^d}{u}^2
      \right)^{\frac12}
      \\
      &\qquad \times
      \left(\int_0^{\tF}
      \seminorm{\beta,h}{\underline{e}_h}^2
      + \nu \norm{1,h}{\underline{e}_h}^2
      \right)^{\frac12}.
    \end{aligned}
  \end{equation}
  Recalling the definition \eqref{eq:def:tnorm} of $\tnorm{h}{{\cdot}}$ we have
  \[
  \begin{gathered}
    \left(\int_0^{\tF} \norm{0,h}{\underline{e}_h}^2\right)^{\frac12}
    \le \tF^{\frac12} \tnorm{h}{\underline{e}_h}\,,\quad
    \left(\int_0^{\tF} \nu \norm{1,h}{\underline{e}_h}^2\right)^{\frac12}\le \tnorm{h}{\underline{e}_h},
    \\
    \left(\int_0^{\tF}
    \seminorm{\beta,h}{\underline{e}_h}^2
    + \nu \norm{1,h}{\underline{e}_h}^2\right)^{\frac12}\le \tnorm{h}{\underline{e}_h}.
  \end{gathered}
  \]
  Plugging these bounds into \eqref{eq:rate.1} and simplifying by $\tnorm{h}{\underline{e}_h}$ concludes the proof. \end{proof}

\begin{remark}[A practical choice for $\beta_T$]
  A possible choice of the stabilization parameter is
  \begin{equation}\label{eq:beta:choice}
    \beta_T = \max\left\{
    c_{\rm s} , \norm{L^\infty(T)^d}{u_T}
    \right\} \qquad \forall T \in {\cal T}_h
  \end{equation}
  where $c_{\rm s} > 0$ is a small ``safeguard'' constant.
  We here show briefly that such choice leads to a uniformly bounded error constant in Corollary~\ref{cor:convergence.rate}, i.e., that the terms depending on $\beta_T$ in the right-hand side of \eqref{eq:coroll.rate} are uniformly bounded.
  Since, by assumption, $\norm{L^\infty(T)^d}{u}$ is bounded for (almost) all time instants, recalling the definition \eqref{eq:chi} of $\chi$ and the choice \eqref{eq:beta:choice} for $\beta_T$, we have
  \[
  \max\left\{ 1, \norm{C^0([0,\tF])}{\chi}\right\} \le
  \max\left\{ 1, \frac{\norm{L^\infty(T)^d}{u}}{c_{\rm s}} \right\} \lesssim 1.
  \]
  Therefore, we are only left to check the term
  \begin{equation}\label{eq:tempL:1}
    \begin{aligned}
      &\sum_{T \in \Th} h_T^{2k+1} \int_0^{\tF} \min(1,\ReT) \, \beta_T \, \seminorm{H^{k+1}(T)^d}{u}^2
      \\
      &\quad
      \lesssim \int_0^{\tF} \sum_{T \in \Th} h_T^{2k+1} \min(1,\ReT) |T|^{\frac12} \beta_T \, \seminorm{W^{k+1,4}(T)^d}{u}^2
      \\
      &\quad
      \lesssim h^{2k+1} \int_0^{\tF} \sum_{T \in \Th} |T|^{\frac12} \beta_T \, \seminorm{W^{k+1,4}(T)^d}{u}^2,
    \end{aligned}
  \end{equation}
  which we bounded by trivial manipulations and the H\"older inequality, assuming some (local) additional regularity for $u$.
  Notice that, in \eqref{eq:tempL:1}, the sum of elements $T$ such that $\beta_T = c_{\rm s}$ is immediate to bound assuming $u \in L^2(0,\tF;W^{k+1,4}({\cal T}_h))$.
  We therefore focus on the elements $T$ such that $\beta_T = \norm{L^\infty(T)^d}{u_T}$ in \eqref{eq:tempL:1} and write
  \[
  \begin{aligned}
    &h^{2k+1}  \int_0^{\tF} \sum_{T \in \Th} |T|^{\frac12} \norm{L^\infty(T)^d}{u_T} \seminorm{W^{k+1,4}(T)^d}{u}^2
    \\
    &\quad
    \lesssim h^{2k+1}  \int_0^{\tF} \sum_{T \in \Th} \norm{L^2(T)^d}{u_T} \seminorm{W^{k+1,4}(T)^d}{u}^2
    \\
    &\quad
    \lesssim h^{2k+1} \int_0^{\tF} \norm{0,h}{\underline{u}_h} \seminorm{W^{k+1,4}({\cal T}_h)^d}{u}^2,
  \end{aligned}
  \]
  where we have used the inverse Lebesgue estimate $\norm{L^\infty(T)}{u_T} \lesssim |T|^{-\frac12} \norm{L^2(T)}{u_T}$ valid for polynomial functions (see, e.g., \cite[Lemma~1.25]{Di-Pietro.Droniou:20})  followed by a discrete H\"older inequality.
  Since Theorem \ref{thm:well.posedness.scheme} implies that $\norm{0,h}{\underline{u}_h}$ is bounded in $L^\infty(0,\tF)$, by again assuming $u \in L^2(0,\tF;W^{k+1,4}({\cal T}_h))$ we can confirm the expected bound.
\end{remark}


\subsection{Estimates of the consistency error components}

\begin{lemma}[Estimate of the unsteady error component]\label{lem:estimate.Err.time}
  Let the assumptions of Theorem~\ref{thm:error.estimate} hold true, and further assume that $u\in H^1(0,\tF;H^{k+1}(\Th)^d)$.
  Then, for almost every $t \in (0,\tF)$ and all $\underline{v}_h \in \UhZ$, it holds, with hidden constant as in Corollary~\ref{cor:convergence.rate} and $\Err_{{\rm time},h}$ defined by \eqref{eq:E.time.h},
  \[ 
  \Err_{{\rm time},h}(\underline{v}_h)
  \lesssim
  \left(
  \sum_{T \in \Th} h_T^{2(k+1)}
  \Seminorm{H^{k+1}(T)^d}{\frac{d u}{dt}}^2
  \right)^{\frac12} \norm{0,h}{\underline{v}_h}.
  \] 
\end{lemma}

\begin{proof}
  Expanding $(\cdot,\cdot)_{0,h}$ according to its definition \eqref{eq:0.h.product}, we get
  \[
  \begin{aligned}
    &\Err_{{\rm time},h}(\underline{v}_h)
    \\
    &\quad
    = \sum_{T \in \Th} \left[
      \int_T \frac{d (u - \widehat{u}_T)}{dt} \cdot v_T
      - h_T \sum_{F \in \FT} \int_F \frac{d(\widehat{u}_F - \widehat{u}_T)}{dt} \cdot (v_F - v_T)
      \right]
    \\
    &\quad
    \le
    \left[
      \sum_{T \in \Th} \left(
      \Norm{L^2(T)^d}{\frac{d( u - \widehat{u}_T )}{dt}}^2
      + h_T \sum_{F \in \FT} \Norm{L^2(F)^d}{\frac{d ( \widehat{u}_F - \widehat{u}_T )}{dt}}^2
      \right)
      \right]^{\frac12} \norm{0,h}{\underline{v}_h}
    \\
    &\quad
    \lesssim
    \left(
    \sum_{T \in \Th} h_T^{2(k+1)}
    \Seminorm{H^{k+1}(T)^d}{\frac{d u}{dt}}^2
    \right)^{\frac12} \norm{0,h}{\underline{v}_h},
  \end{aligned}
  \]
  where we have used Cauchy--Schwarz inequalities in the second passage.
  To obtain the conclusion, for the first term we have used the approximation properties \eqref{eq:IRTNT:approximation} of the Raviart--Thomas--Nédélec interpolator with $(s,\ell,q,m) = (2,k+1,k,0)$ after noticing that $\frac{d \widehat{u}_T}{dt} = \IRTNT{k+1} \left(\frac{d u}{dt}\right)$;
  for the second term, we proceeded in a similar way as for \eqref{eq:W.1.infty.boundedness.IUh:T2}, except for the use of $L^2$-norms instead of $L^\infty$-norms and the fact that the higher regularity of $\frac{d u}{dt}$ is exploited when invoking the approximation properties of $\IRTNT{k+1}$ and $\lproj{k}{T}$.
\end{proof}

The following bound for the diffusive error is proved in \cite[Theorem~14]{Botti.Botti.ea:24} using standard arguments from the Hybrid High-Order literature; see, e.g., \cite[Point (ii) in Lemma 2.18]{Di-Pietro.Droniou:20}.

\begin{lemma}[Estimate of the diffusive error component]\label{lem:estimate.Err.diff.h}
  Let the assumptions of Theorem~\ref{thm:error.estimate} hold true, and further assume that $u \in L^2(0,\tF;H^{k+2}(\Th))^d$.
  Then, for almost every $t \in (0,\tF)$ and all $\underline{v}_h \in \UhZ$, it holds, with hidden constant as in Corollary~\ref{cor:convergence.rate} and $\Err_{{\rm diff},h}$ defined by \eqref{eq:E.diff.h},
  \[ 
  \Err_{{\rm diff},h}(\underline{v}_h)
  \lesssim
  \left(
  \sum_{T \in \Th} \nu h_T^{2(k+1)} \seminorm{H^{k+2}(T)^d}{u}^2
  \right)^{\frac12}
  \nu^{\frac12} \norm{1,h}{\underline{v}_h}.
  \] 
\end{lemma}

\begin{lemma}[Regime-dependent estimate of the convective error component]\label{lem:Err.conv.h:estimate}
  Let the assumptions of Theorem~\ref{thm:error.estimate} hold true, and further assume that $u \in L^2(0,\tF;H^{k+1}(\Th))^d$ and that the quantity defined in \eqref{eq:chi} satisfies $\chi \in C^0([0,\tF])$.
  Then, for almost every $t \in (0,\tF)$ and all $\underline{v}_h \in \ZhZ$, it holds, with hidden constant as in Corollary~\ref{cor:convergence.rate}  and $\Err_{{\rm conv},h}$ defined by \eqref{eq:E.conv.h},
  \begin{multline}\label{eq:estimate.Err.conv.h}
    \Err_{{\rm conv},h}(\underline{v}_h)
    \lesssim
    \left(
    \sum_{T \in \Th} h_T^{2(k+1)}
    \seminorm{W^{1,\infty}(T)^d}{u}^2 \seminorm{H^{k+1}(T)^d}{u}^2
    \right)^{\frac12} \norm{L^2(\Omega)^d}{v_h}
    \\
    +  ( 1 + \norm{C^0([0,\tF])}{\chi})^{\frac12}
    \left(
    \sum_{T \in \Th} h_T^{2k+1} \min(1,\ReT)
    \big(
    \beta_T +  \norm{L^\infty(T)^d}{u}
    \big) \seminorm{H^{k+1}(T)^d}{u}^2
    \right)^{\frac12}
    \\
    \times \left(
    \seminorm{\beta,h}{\underline{v}_h}^2
    + \nu \norm{1,h}{\underline{v}_h}^2
    \right)^{\frac12}.
  \end{multline}
\end{lemma}

\begin{proof}
  Expanding $t_h$ according to its definition \eqref{eq:th} and writing $v_F+v_T=2v_T+(v_F-v_T)$ in the boundary term, we decompose the convective error as follows:
  \begin{equation}\label{eq:Err.conv:estimate:basic}
    \begin{aligned}
      &\Err_{{\rm conv},h}(\underline{v}_h)
      \\
      &\quad =
      \sum_{T \in \Th} \left(
      \int_T (u \cdot \nabla) u \cdot v_T
      - \int_T (\widehat{u}_T \cdot \nabla) \widehat{u}_T \cdot v_T
      - \sum_{F \in \FT} \!\! \int_F (\widehat{u}_T \cdot n_{TF}) (\widehat{u}_F - \widehat{u}_T) \cdot v_T
      \right)
      \\
      &\qquad
      - \sum_{T \in \Th} \sum_{F \in \FT} \int_F
      \left[
        \beta_T + \frac12(\widehat{u}_T \cdot n_{TF})
        \right]
      (\widehat{u}_F - \widehat{u}_T) \cdot (v_F - v_T)
      \\
      &\quad \eqcolon
      \term_1 + \term_2.
    \end{aligned}
  \end{equation}

  In what follows, we write $\term_i = \sum_{T \in \Th} \term_i(T)$ and estimate the error contribution for a given mesh element $T \in \Th$.
  To estimate $\term_{1}(T)$, we let $\overline{u}_T \coloneq \lproj{0}{T} \widehat{u}_T$ and write
  \[
  \begin{aligned}
    \int_T (\overline{u}_T \cdot \nabla) u \cdot v_T
    &=
    - \int_T u \cdot \nabla \cdot( v_T \otimes \overline{u}_T)
    + \sum_{F \in \FT} \int_F (\overline{u}_T \cdot n_{TF}) (u \cdot v_T)
    \\
    \overset{\eqref{eq:lproj}}&=
    - \int_T \lproj{k-1}{T} u \cdot \nabla \cdot(  v_T \otimes \overline{u}_T)
    + \sum_{F \in \FT} \int_F (\overline{u}_T \cdot n_{TF}) (\lproj{k}{F} u \cdot  v_T)
    \\
    &=
    - \int_T \widehat{u}_T \cdot \nabla \cdot(  v_T \otimes \overline{u}_T)
    + \sum_{F \in \FT} \int_F (\overline{u}_T \cdot n_{TF}) (\widehat{u}_F \cdot  v_T),
  \end{aligned}
  \]
  where we have used an integration by parts in the first equality, we have noticed that ${\nabla \cdot (v_T \otimes \overline{u}_T)} \in \Poly{k-1}(T)^d$ (by \eqref{eq:ZhZ.degree.T}) and $v_T {(\overline{u}_T \cdot n_{TF})} \in \Poly{k}(F)^d$ to insert the projectors in the second equality,
  while, in the third equality, we have first written
  $\lproj{k-1}{T} u \overset{\eqref{eq:IRTNT},\,\eqref{eq:IUh}} = \lproj{k-1}{T} \widehat{u}_T$
  and then removed $\lproj{k-1}{T}$ using again ${\nabla \cdot (v_T \otimes \overline{u}_T)} \in \Poly{k-1}(T)^d$.
  Integrating by parts the first term in right-hand side of the above expression and rearranging, we obtain
  \[
  \int_T ( \overline{u}_T \cdot \nabla ) (u - \widehat{u}_T) \cdot v_T
  - \sum_{F \in \FT} \int_F (\overline{u}_T \cdot n_{TF}) (\widehat{u}_F - \widehat{u}_T) \cdot v_T = 0.
  \]
  Subtracting the above quantity from $\term_{1}(T)$
  and adding and subtracting $\int_T (\widehat{u}_T \cdot \nabla) u \cdot v_T$ to the resulting expression, we obtain
  \[
  \begin{aligned}
    &\term_{1}(T)
    \\
    &\quad
    = \int_T [ (u - \widehat{u}_T) \cdot \nabla] u \cdot v_T
    + \int_T [ (\widehat{u}_T - \overline{u}_T) \cdot \nabla ] (u - \widehat{u}_T) \cdot v_T
    \\
    &\qquad
    - \sum_{F \in\ FT} \int_F (\widehat{u}_T - \overline{u}_T) \cdot n_{TF} (\widehat{u}_F - \widehat{u}_T) \cdot v_T
    \\
    &\quad\lesssim
    \left(
    \seminorm{W^{1,\infty}(T)^d}{u} \norm{L^2(T)^d}{u {-} \widehat{u}_T}
    + \norm{L^2(T)^{d \times d}}{\nabla (u {-} \widehat{u}_T)}
    \norm{L^\infty(T)^d}{\widehat{u}_T {-} \overline{u}_T}
    \right) \norm{L^2(T)^d}{v_T}
    \\
    &\qquad
    + \sum_{F \in \FT}
    h_T^{-1}
    \norm{L^\infty(F)^d}{\widehat{u}_T - \overline{u}_T}
    ~h_T^{\frac12} \norm{L^2(F)^d}{\widehat{u}_F - \widehat{u}_T}
    ~h_T^{\frac12} \norm{L^2(F)^d}{v_T},
  \end{aligned}
  \]
  where we have used H\"{o}lder inequalities in the second step.
  We next proceed to bound the various norms of approximation errors on $u$ that appear in the above expression.
  Using \eqref{eq:IRTNT:approximation} with ${(s,\ell,q,m)} = {(2,k+1,k,0)}$, we have
  \[
  \norm{L^2(T)^d}{u - \widehat{u}_T} \lesssim h_T^{k+1} \seminorm{H^{k+1}(T)^d}{u}.
  \]
  Recalling that $\overline{u}_T=\lproj{0}{T}\widehat{u}_T$ and using the approximation properties of $\lproj{0}{T}$ we write
  \begin{equation*}
    {\norm{L^\infty(T)^d}{\widehat{u}_T - \overline{u}_T}}
    \lesssim h_T \seminorm{W^{1,\infty}(T)^d}{\widehat{u}_T}
    \overset{\eqref{eq:boundedness.IRTNT}}\lesssim h_T \seminorm{W^{1,\infty}(T)^d}{u}.
  \end{equation*}
  By Proposition \ref{prop:Interpolate.ZhZ}, we have $\underline{u}_h \in \ZhZ$ and thus, applying \eqref{eq:ZhZ.degree.T}, $\widehat{u}_T\in\Poly{k}(T)^d$. By definition of $\widehat{u}_F = \lproj{k}{F}u$ and continuity of $\lproj{k}{F}$, we infer
  \begin{equation}\label{eq:estimate.uF-pi.k.T.hu.T}
    \norm{L^2(F)^d}{\widehat{u}_F - \widehat{u}_T} =
    \norm{L^2(F)^d}{\lproj{k}{F}(u - \widehat{u}_T)}
    \le \norm{L^2(F)^d}{u - \widehat{u}_T}
    \overset{\eqref{eq:IRTNT:approximation.trace}}\lesssim h_T^{k+\frac12} \seminorm{H^{k+1}(T)^d}{u} .
  \end{equation}
  By collecting all of the above bounds and observing that, by a discrete trace inequality
  $h_T^{\frac12} {\norm{L^2(F)^d}{v_T}} \lesssim {\norm{L^2(T)^d}{v_T}} $, we eventually find
  \begin{equation*}
    \term_{1}(T) \lesssim
    h_T^{k+1} \seminorm{W^{1,\infty}(T)^d}{u}
    \seminorm{H^{k+1}(T)^d}{u}
    \norm{L^2(T)^d}{v_T}.
  \end{equation*}
  Using a Cauchy--Schwarz inequality on the sum over the elements, we conclude that
  \begin{equation}\label{eq:Err.conv:estimate:T1}
    \term_1
    \lesssim \left[
      \sum_{T \in \Th} h_T^{2(k+1)}
      \seminorm{W^{1,\infty}(T)^d}{u}^2 \seminorm{H^{k+1}(T)^d}{u}^2
      \right]^{\frac12} \norm{L^2(\Omega)^d}{v_h}.
  \end{equation}

  To estimate $\term_2(T)$, we start with Cauchy--Schwarz and $(\infty,2,2)$-H\"{o}lder inequalities to infer
  \[
  \term_2(T)
  \lesssim  \sum_{F \in \FT} \big(
  \beta_T + \norm{L^\infty(F)^d}{\widehat{u}_T \cdot n_{TF}}
  \big) \norm{L^2(F)^d}{\widehat{u}_F - \widehat{u}_T}
  \norm{L^2(F)^d}{v_F - v_T}.
  \]
  We next notice that
  \[
    {\norm{L^\infty(F)}{\widehat{u}_T \cdot n_{TF}}}
    \overset{\eqref{eq:IUh},\eqref{eq:IRTNT}}= {\norm{L^\infty(F)}{\lproj{k}{F}(u \cdot n_{TF})}}
    \lesssim {\norm{L^\infty(F)}{u \cdot n_{TF}}}
    \le {\norm{L^\infty(T)^d}{u}},
    \]
    where we have used the $L^\infty$-boundedness of $\lproj{k}{F}$ in the first inequality.
    Accounting for this fact, recalling \eqref{eq:estimate.uF-pi.k.T.hu.T}, and using a Cauchy--Schwarz inequality on the sum over the faces, we get
    \begin{multline*}
      \term_2(T)
      \lesssim h_T^{k+\frac12}
      \big(
      \beta_T + \norm{L^\infty(T)^d}{u}
      \big)^{\frac12} \seminorm{H^{k+1}(T)^d}{u}
      \\
      \times
      \left[
        \big(
        \beta_T + \norm{L^\infty(T)^d}{u}
        \big) \sum_{F \in \FT}\norm{L^2(F)^d}{v_F - v_T}^2
        \right]^{\frac12}.
    \end{multline*}
    We next proceed differently according to the regime.
    If $\ReT > 1$, then we recall the definitions \eqref{eq:seminorm.beta} of the $\seminorm{\beta,T}{{\cdot}}$-seminorm and \eqref{eq:chi} of $\chi$ to arrive at
    \begin{equation}\label{eq:estimate.Err.conv:T2(T).conv}
      \term_2(T)
      \lesssim h_T^{k+\frac12}
      \big(
      \beta_T + \norm{L^\infty(T)^d}{u}
      \big)^{\frac12} \seminorm{H^{k+1}(T)^d}{u}
      (1 + \chi)^{\frac12} \seminorm{\beta,T}{\underline{v}_T}
      \qquad \forall T \in \Th,\,\ReT > 1.
    \end{equation}
    If, on the other hand, $\ReT \le 1$, we start by noticing that
    \begin{multline*}
    \big(
    \beta_T + \norm{L^\infty(T)^d}{u}
    \big) \sum_{F \in \FT}\norm{L^2(F)^d}{v_F - v_T}^2
    \\
    \overset{\eqref{eq:ReT}}= \ReT \nu h_T^{-1} \sum_{F \in \FT}\norm{L^2(F)^d}{v_F - v_T}^2
    \le \nu \ReT \norm{1,T}{\underline{v}_T}^2,
    \end{multline*}
    where the conclusion follows from the definition \eqref{eq:norm.1.T} of $\norm{1,T}{{\cdot}}$.
    Hence,
    \begin{equation}\label{eq:estimate.Err.conv:T2(T).diff}
      \term_2(T) \lesssim h_T^{k+\frac12} \ReT^{\frac12}
      \big(
      \beta_T + \norm{L^\infty(T)^d}{u}
      \big)^{\frac12} \seminorm{H^{k+1}(T)^d}{u}~\nu^{\frac12} \norm{1,T}{\underline{v}_T}\qquad \forall T \in \Th,\,\ReT \le 1.
    \end{equation}
    Gathering the estimates \eqref{eq:estimate.Err.conv:T2(T).conv} and \eqref{eq:estimate.Err.conv:T2(T).diff} leads to
    \[
    \term_2(T) \lesssim (1 + \chi)h_T^{k+\frac12} \!
    \min(1,\ReT)^{\frac12}\big(\beta_T + \norm{L^\infty(T)^d}{u}\big)^{\frac12}\seminorm{H^{k+1}(T)^d}{u} \!
    \left( \seminorm{\beta,T}{\underline{v}_T}^2 {+} \nu\norm{1,T}{\underline{v}_T}^2\right)^{\frac12}
    \]
    for all $T\in\Th$.
    Using a Cauchy--Schwarz inequality on the sum, we arrive at
    \begin{multline}\label{eq:Err.conv:estimate:T2}
      \term_2
      \lesssim
      ( 1 + \norm{C^0([0,\tF])}{\chi})^{\frac12}
      \left(
      \sum_{T \in \Th} h_T^{2k+1} \min(1, \ReT)
      \big(
      \beta_T + \norm{L^\infty(T)^d}{u}
      \big) \seminorm{H^{k+1}(T)^d}{u}^2
      \right)^{\frac12}
      \\
      \times \left(
      \seminorm{\beta,h}{\underline{v}_h}^2
      + \nu \norm{1,h}{\underline{v}_h}^2
      \right)^{\frac12}.
    \end{multline}
    The conclusion follows using \eqref{eq:Err.conv:estimate:T1} and \eqref{eq:Err.conv:estimate:T2} to estimate the right-hand side of \eqref{eq:Err.conv:estimate:basic}.
\end{proof}



\section{Numerical tests}\label{sec:num}

In this section, we numerically verify the proposed method applied to a two-dimensional analytic test problem studied in \cite{Castanon-Quiroz.Di-Pietro:25,Han.Hou:21,de-frutos-fully-2019} on a spatial domain $\Omega=(0,1)^2$ and final time $\tF=1$. The forcing data $f$ is chosen such that the velocity and pressure components are such that
\begin{align*}
  u(t,x,y) &= \frac{3+2\cos(4t)}{5} \begin{pmatrix}
    16y(1-y)(1-2y)\sin^2(\pi x) \\
    -8\pi y^2(1-y)^2\sin(2\pi x)
  \end{pmatrix}, \\
  p(t,x,y) &= \frac{3+2\cos(4t)}{5} \sin(\pi x) \cos(\pi y).
\end{align*}
We test the method on a family of $h$-refined triangular meshes and a range of polynomial degrees $k\in\{0,1,2\}$. We employ the choice of $\beta_T$ in \eqref{eq:beta:choice} with $c_s=10^{-4}$ and the values of $u_T$ taken at the current time-step in the convective stabilization. An implicit Crank--Nicolson temporal discretization is employed, using the interpolate of the analytic solution as the initial condition. The resulting nonlinear problem for each time-step is solved with the Newton algorithm with a solver tolerance of $10^{-8}$. We notice that, since the incompressibility equation \eqref{eq:discrete:mass} is linear, the divergence-free condition is exactly satisfied at \emph{each step} of the Newton iteration, and is therefore not impacted by the choice of the tolerance.
The total number of time-steps is taken to be
\begin{equation*}
  N_{\tF} = \max\left\{ 10, \left\lceil \frac{1}{h^{(k+1)/2}} \right\rceil \right\},
\end{equation*}
where $h=\max_{T\in\Th}h_T$; this choice ensures that the global time-integration error using the time-step size $\Delta t = \frac{\tF}{N_{\tF}}$ is of order $h^{k+1}$. To demonstrate the robustness of the velocity error estimate with respect to the local Reynolds number $\ReT$ defined in \eqref{eq:ReT}, we test a range of values for $\nu\in\{10^{-6},10^{-2},1\}$, with $\nu=1$ corresponding to a diffusion-dominated ($\ReT<1$) regime and $\nu=10^{-6}$ corresponding to a convection-dominated ($\ReT>1$) regime on the sequence of meshes considered.

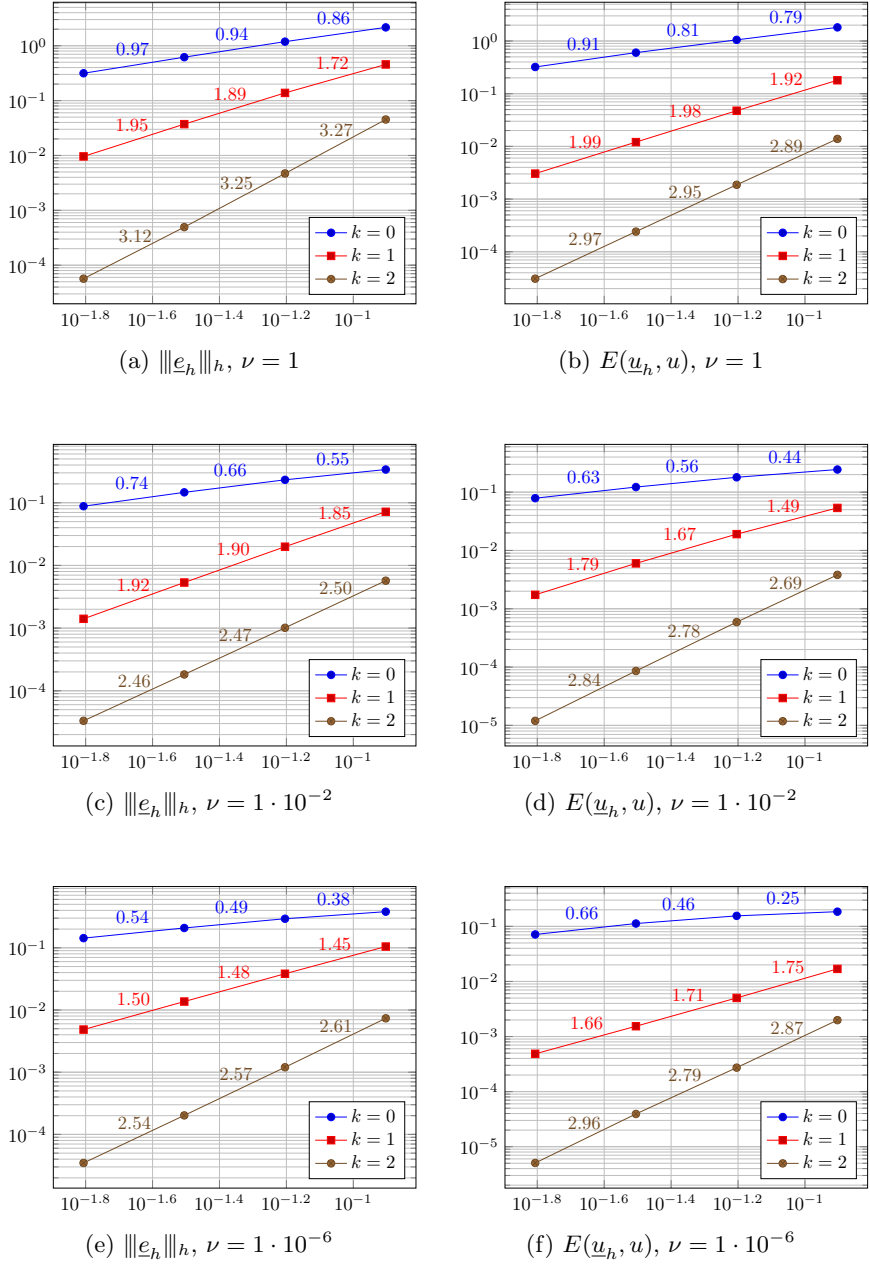
\begin{figure}\centering

  \subcaptionbox{$\tnorm{h}{\underline{e}_h}$, $\nu = \pgfmathprintnumber{1.}$}[0.45\textwidth]{\centering
    \begin{tikzpicture}[scale=0.70]
      \begin{loglogaxis}[legend pos=south east, grid=both]
        \addplot table[col sep=comma,x=MeshSize,y=TotalEnergyComponentNormError] {outputs/tri_k0_mu1./data_rates.dat}
        node[pos = 0.166666666666667, above=2pt]{0.86} 
        node[pos = 0.5, above=2pt]{0.94} 
        node[pos = 0.833333333333333, above=2pt]{0.97}; 
        \addlegendentry{$k=0$}

        \addplot table[col sep=comma,x=MeshSize,y=TotalEnergyComponentNormError] {outputs/tri_k1_mu1./data_rates.dat}
        node[pos = 0.166666666666667, above=2pt]{1.72} 
        node[pos = 0.5, above=2pt]{1.89} 
        node[pos = 0.833333333333333, above=2pt]{1.95}; 
        \addlegendentry{$k=1$}

        \addplot table[col sep=comma,x=MeshSize,y=TotalEnergyComponentNormError] {outputs/tri_k2_mu1./data_rates.dat}
        node[pos = 0.166666666666667, above=2pt]{3.27} 
        node[pos = 0.5, above=2pt]{3.25} 
        node[pos = 0.833333333333333, above=2pt]{3.12}; 
        \addlegendentry{$k=2$}

      \end{loglogaxis}
    \end{tikzpicture}
  }
  \subcaptionbox{$E(\underline{u}_h, u)$, $\nu = \pgfmathprintnumber{1.}$}[0.45\textwidth]{\centering
    \begin{tikzpicture}[scale=0.70]
      \begin{loglogaxis}[legend pos=south east, grid=both]
        \addplot table[col sep=comma,x=MeshSize,y=TotalEnergyError] {outputs/tri_k0_mu1./data_rates.dat}
        node[pos = 0.166666666666667, above=2pt]{0.79} 
        node[pos = 0.5, above=2pt]{0.81} 
        node[pos = 0.833333333333333, above=2pt]{0.91}; 
        \addlegendentry{$k=0$}

        \addplot table[col sep=comma,x=MeshSize,y=TotalEnergyError] {outputs/tri_k1_mu1./data_rates.dat}
        node[pos = 0.166666666666667, above=2pt]{1.92} 
        node[pos = 0.5, above=2pt]{1.98} 
        node[pos = 0.833333333333333, above=2pt]{1.99}; 
        \addlegendentry{$k=1$}

        \addplot table[col sep=comma,x=MeshSize,y=TotalEnergyError] {outputs/tri_k2_mu1./data_rates.dat}
        node[pos = 0.166666666666667, above=2pt]{2.89} 
        node[pos = 0.5, above=2pt]{2.95} 
        node[pos = 0.833333333333333, above=2pt]{2.97}; 
        \addlegendentry{$k=2$}

      \end{loglogaxis}
    \end{tikzpicture}
  }
  \\ \vspace{0.5cm}
  \subcaptionbox{$\tnorm{h}{\underline{e}_h}$, $\nu = \pgfmathprintnumber{1.e-2}$}[0.45\textwidth]{\centering
    \begin{tikzpicture}[scale=0.70]
      \begin{loglogaxis}[legend pos=south east, grid=both]
        \addplot table[col sep=comma,x=MeshSize,y=TotalEnergyComponentNormError] {outputs/tri_k0_mu1.e-2/data_rates.dat}
        node[pos = 0.166666666666667, above=2pt]{0.55} 
        node[pos = 0.5, above=2pt]{0.66} 
        node[pos = 0.833333333333333, above=2pt]{0.74}; 
        \addlegendentry{$k=0$}

        \addplot table[col sep=comma,x=MeshSize,y=TotalEnergyComponentNormError] {outputs/tri_k1_mu1.e-2/data_rates.dat}
        node[pos = 0.166666666666667, above=2pt]{1.85} 
        node[pos = 0.5, above=2pt]{1.90} 
        node[pos = 0.833333333333333, above=2pt]{1.92}; 
        \addlegendentry{$k=1$}

        \addplot table[col sep=comma,x=MeshSize,y=TotalEnergyComponentNormError] {outputs/tri_k2_mu1.e-2/data_rates.dat}
        node[pos = 0.166666666666667, above=2pt]{2.50} 
        node[pos = 0.5, above=2pt]{2.47} 
        node[pos = 0.833333333333333, above=2pt]{2.46}; 
        \addlegendentry{$k=2$}

      \end{loglogaxis}
    \end{tikzpicture}
  }
  \subcaptionbox{$E(\underline{u}_h, u)$, $\nu = \pgfmathprintnumber{1.e-2}$}[0.45\textwidth]{\centering
    \begin{tikzpicture}[scale=0.70]
      \begin{loglogaxis}[legend pos=south east, grid=both]
        \addplot table[col sep=comma,x=MeshSize,y=TotalEnergyError] {outputs/tri_k0_mu1.e-2/data_rates.dat}
        node[pos = 0.166666666666667, above=2pt]{0.44} 
        node[pos = 0.5, above=2pt]{0.56} 
        node[pos = 0.833333333333333, above=2pt]{0.63}; 
        \addlegendentry{$k=0$}

        \addplot table[col sep=comma,x=MeshSize,y=TotalEnergyError] {outputs/tri_k1_mu1.e-2/data_rates.dat}
        node[pos = 0.166666666666667, above=2pt]{1.49} 
        node[pos = 0.5, above=2pt]{1.67} 
        node[pos = 0.833333333333333, above=2pt]{1.79}; 
        \addlegendentry{$k=1$}

        \addplot table[col sep=comma,x=MeshSize,y=TotalEnergyError] {outputs/tri_k2_mu1.e-2/data_rates.dat}
        node[pos = 0.166666666666667, above=2pt]{2.69} 
        node[pos = 0.5, above=2pt]{2.78} 
        node[pos = 0.833333333333333, above=2pt]{2.84}; 
        \addlegendentry{$k=2$}

      \end{loglogaxis}
    \end{tikzpicture}
  }
  \\ \vspace{0.5cm}
  \subcaptionbox{$\tnorm{h}{\underline{e}_h}$, $\nu = \pgfmathprintnumber{1.e-6}$}[0.45\textwidth]{\centering
    \begin{tikzpicture}[scale=0.70]
      \begin{loglogaxis}[legend pos=south east, grid=both]
        \addplot table[col sep=comma,x=MeshSize,y=TotalEnergyComponentNormError] {outputs/tri_k0_mu1.e-6/data_rates.dat}
        node[pos = 0.166666666666667, above=2pt]{0.38} 
        node[pos = 0.5, above=2pt]{0.49} 
        node[pos = 0.833333333333333, above=2pt]{0.54}; 
        \addlegendentry{$k=0$}

        \addplot table[col sep=comma,x=MeshSize,y=TotalEnergyComponentNormError] {outputs/tri_k1_mu1.e-6/data_rates.dat}
        node[pos = 0.166666666666667, above=2pt]{1.45} 
        node[pos = 0.5, above=2pt]{1.48} 
        node[pos = 0.833333333333333, above=2pt]{1.50}; 
        \addlegendentry{$k=1$}

        \addplot table[col sep=comma,x=MeshSize,y=TotalEnergyComponentNormError] {outputs/tri_k2_mu1.e-6/data_rates.dat}
        node[pos = 0.166666666666667, above=2pt]{2.61} 
        node[pos = 0.5, above=2pt]{2.57} 
        node[pos = 0.833333333333333, above=2pt]{2.54}; 
        \addlegendentry{$k=2$}

      \end{loglogaxis}
    \end{tikzpicture}
  }
  \subcaptionbox{$E(\underline{u}_h, u)$, $\nu = \pgfmathprintnumber{1.e-6}$}[0.45\textwidth]{\centering
    \begin{tikzpicture}[scale=0.70]
      \begin{loglogaxis}[legend pos=south east, grid=both]
        \addplot table[col sep=comma,x=MeshSize,y=TotalEnergyError] {outputs/tri_k0_mu1.e-6/data_rates.dat}
        node[pos = 0.166666666666667, above=2pt]{0.25} 
        node[pos = 0.5, above=2pt]{0.46} 
        node[pos = 0.833333333333333, above=2pt]{0.66}; 
        \addlegendentry{$k=0$}

        \addplot table[col sep=comma,x=MeshSize,y=TotalEnergyError] {outputs/tri_k1_mu1.e-6/data_rates.dat}
        node[pos = 0.166666666666667, above=2pt]{1.75} 
        node[pos = 0.5, above=2pt]{1.71} 
        node[pos = 0.833333333333333, above=2pt]{1.66}; 
        \addlegendentry{$k=1$}

        \addplot table[col sep=comma,x=MeshSize,y=TotalEnergyError] {outputs/tri_k2_mu1.e-6/data_rates.dat}
        node[pos = 0.166666666666667, above=2pt]{2.87} 
        node[pos = 0.5, above=2pt]{2.79} 
        node[pos = 0.833333333333333, above=2pt]{2.96}; 
        \addlegendentry{$k=2$}

      \end{loglogaxis}
    \end{tikzpicture}
  }
  \caption{Convergence rates for the numerical tests of Section \ref{sec:num} with $\nu\in\{10^{-6},10^{-2},1\}$ and $k\in\{0,1,2\}$. Estimated orders of convergence between successive spatial $h$-refinements are also given.\label{fig:convergence:tnorm}}
\end{figure}

The left column of Figure \ref{fig:convergence:tnorm} shows the behaviour of the velocity error $\underline{e}_h$ defined  in \eqref{eq:eh} measured in a time-discrete version of the norm $\tnorm{h}{{\cdot}}$ in \eqref{eq:def:tnorm}. The estimated orders of convergence of the method are also presented in Figure \ref{fig:convergence:tnorm}. Recalling the statement of Corollary \ref{cor:convergence.rate}, we expect to observe an error reduction rate of $k+\frac12$ in the convection-dominated case and $k+1$ in the diffusion-dominated case, which is indeed shown in Figure \ref{fig:convergence:tnorm}. The transition from pre-asymptotic to asymptotic error reduction rates as the mesh is refined (and thus as $\ReT$ decreases) is also observed for the intermediate value $\nu = 10^{-2}$.
For the sake of completeness, we also report in the right column of Figure \ref{fig:convergence:tnorm} a time-discrete version of the error with respect to the continuous solution:
\[
E(\underline{u}_h, u)
\coloneq
\left(
\max_{t \in [0, \tF]} \norm{L^2(\Omega)^d}{\Rh \underline{u}_h - u}^2
+ \int_0^{\tF} \nu \norm{L^2(\Omega)^{d\times d}}{\nabla_h \Rh \underline{u}_h - \nabla u}^2
\right)^{\frac12},
\]
where the global velocity reconstruction $\Rh \underline{u}_h$ is such that $(\Rh \underline{u}_h)_{|T} \coloneq \RT \underline{u}_T$ for all $T \in \Th$ and $\nabla_h$ denotes the broken gradient on $\Th$.
This error can be bounded using triangle inequalities, discrete norm equivalences, \eqref{eq:coroll.rate}, and the approximation properties of the velocity reconstruction.
Overall, $E(\underline{u}_h, u)$ displays the same behaviour as $\tnorm{h}{\underline{e}_h}$, with slightly higher convergence rates in some cases, which can be attributed to the faster convergence of the approximation component.

\section*{Acknowledgements}

Funded by the European Union (ERC Synergy, NEMESIS, project number 101115663).
Views and opinions expressed are however those of the authors only and do not necessarily reflect those of the European Union or the European Research Council Executive Agency. Neither the European Union nor the granting authority can be held responsible for them.


\printbibliography

\end{document}